\documentclass{article}%
\usepackage{graphicx}
\usepackage{amsfonts}
\usepackage{amssymb}
\usepackage{amsmath}%
\providecommand{\U}[1]{\protect\rule{.1in}{.1in}}
\begin{document}

\title{Measuring and Controlling Bias for Some Bayesian Inferences and the Relation
to Frequentist Criteria}
\author{Michael Evans and Yang Guo\\Department of Statistical Sciences\\University of Toronto}
\date{}
\maketitle

\begin{center}
\textbf{Abstract}
\end{center}

\noindent A common concern with Bayesian methodology in scientific contexts is
that inferences can be heavily influenced by subjective biases. As presented
here, there are two types of bias for some quantity of interest: bias against
and bias in favor. Based upon the principle of evidence, it is shown how to
measure and control these biases for both hypothesis assessment and estimation
problems. Optimality results are established for the principle of evidence as
the basis of the approach to these problems. A close relationship is
established between measuring bias in Bayesian inferences and frequentist
properties that hold for any proper prior. This leads to a possible resolution
to an apparent conflict between these approaches to statistical reasoning.
Frequentism is seen as establishing a figure of merit for a statistical study,
while Bayesianism plays the key role in determining inferences based upon
statistical evidence.\medskip


\section{Introduction}

A serious concern with Bayesian methodology is that the choice of the prior
could result in conclusions that to some degree are predetermined before
seeing the data. In certain circumstances this is correct. This can be seen by
considering the problem associated with what is known as the Jeffreys-Lindley
paradox where posterior probabilities of hypotheses, as well as associated
Bayes factors, will produce increasing support for the hypothesis as the prior
becomes more diffuse. So, while one may feel that a very diffuse prior is
putting in very little information, it is in fact biasing the results in favor
of the hypothesis. It has been argued, see Baskurt and Evans (2013) and Evans
(2015), that the measurement and control of bias is a key element of a
Bayesian analysis as without it, and the assurance that bias is minimal, the
validity of any inference is suspect.

While attempts have been made to avoid the Jeffreys-Lindley paradox through
the choice of the prior, modifying the prior to avoid bias is contrary to the
ideals of a Bayesian analysis which requires the elicitation of a prior based
upon knowledge of the phenomenon under study. Why should one change such a
prior because of bias? Indeed, as will be discussed, there is bias in favor
and bias against and typically choosing a prior to minimize one type of bias
simply increases the other. The real method for controlling bias of both types
is through the amount of data collected. So controlling bias is an aspect of
design. Bias can be measured post-hoc and it then provides a way to assess the
weight that should be given the results of an analysis. For example, if a
study concludes that there is evidence in favor of a hypothesis, but it can be
shown that there was a high prior probability that such evidence would be
obtained, then the results of such an analysis can't be considered to be reliable.

Previous discussion concerning bias was focused on hypothesis assessment and
in many ways this is a natural starting point. This paper is concerned with
adding some aspects to those developments and to extending the approach to
estimation and prediction problems. Furthermore, it is shown here that
measuring and controlling bias establishes close links between a frequentist
approach to statistics and Bayesian inference. In essence frequentism is
concerned with design while inferences are Bayesian. Bayesian inference is
based upon the evidence in the observed data and is unconcerned, at least for
inference, about data sets that could have been obtained. Frequentism is
concerned with the behavior of inferences as applied to unobserved data sets
and this is entirely appropriate before the data is observed. So consideration
of bias leads to a degree of unification between different ways of thinking
about statistical reasoning.

The measurement of bias, and thus its control, is dependent upon measuring
evidence. The \textit{principle of evidence} is adopted here: evidence in
favor of a specific value of an unknown occurs when the posterior probability
of the value is greater than its prior probability, evidence against occurs
when the posterior probability of the value is less than its prior probability
and there is no evidence either way when these are equal. The major part of
what is discussed here depends only on this simple principle but sometimes a
numerical measure of evidence is needed and for this we use the
\textit{relative belief ratio} defined as the ratio of the posterior to prior
probability. The relative belief ratio is related to the Bayes factor but has
some nicer properties such as providing a measure of the evidence for each
value of a parameter without the need to modify the prior.

There is not much discussion in the Bayesian literature of the notion of bias
in the sense that is meant here. There is considerable discussion, however,
concerning the Jeffreys-Lindley paradox and our position is that bias plays a
key role in the issues that arise. Relevant recent papers on this include
Shafer (1982), Spanos (2013), Sprenger (2013), Robert (2014), Cousins (2017)
and Villa and Walker (2017) and these contain extensive background references.
Gu et al. (2019) is concerned with the validation of quantum theory using
Bayesian methodology applied to well-known data sets\ and the principle of
evidence and an assessment of the bias in the prior plays a key role in the argument.

In Section 2 the concepts are defined, their properties are considered and
illustrated via a simple example where the Jeffreys-Lindley paradox is
relevant. Also, it is seen that a well-known p-value does not satisfy the
principle of evidence but can still be used to characterize evidence for or
against but requires significance levels that go to 0 with increasing sample
size or increasing diffuseness of the prior. In Section 3 the relationship
with frequentism is discussed and a number of optimality results are
established for the approach taken here to measuring and controlling bias,
namely, via the principle of evidence. In Section 4, a variety of examples are
considered and analyzed from the point-of-view of bias. All proofs of theorems
are in the Appendix.

\section{Evidence and Bias}

For the discussion here there is a model $\{f_{\theta}:\theta\in\Theta\},$
given by densities $f_{\theta},$ for data $x$ and a proper prior probability
distribution given by density $\pi.$ It is supposed that interest is in
inferences about $\psi=\Psi(\theta)$ where $\Psi:\Theta\rightarrow\Psi$ is
onto and for economy the same notation is used for the function and its range.
For the most part it is safe to assume all the probability distributions are
discrete with results for the continuous case obtained by taking limits.

A measure of the evidence that $\psi\in\Psi$ is the true value is given by the
relative belief ratio
\begin{equation}
RB_{\Psi}(\psi\,|\,x)=\lim_{\delta\rightarrow0}\frac{\Pi_{\Psi}(N_{\delta
}(\psi)\,|\,x)}{\Pi_{\Psi}(N_{\delta}(\psi))}=\frac{\pi_{\Psi}(\psi
\,|\,x)}{\pi_{\Psi}(\psi)} \label{eq1}%
\end{equation}
where $\Pi_{\Psi},\Pi_{\Psi}(\cdot\,|\,x)$ are the prior and posterior
probability measures of $\Psi$ with densities $\pi_{\Psi}$ and $\pi_{\Psi
}(\cdot\,|\,x),$ respectively, and $N_{\delta}(\psi)$ is a sequence of sets
converging nicely to $\{\psi\}.$ The last equality in (\ref{eq1}) requires
some conditions but the prior density positive and continuous at $\psi\,\ $is
enough. So $RB_{\Psi}(\psi\,|\,x)>1$ implies evidence for the true value being
$\psi,$ etc. Any \textit{valid} measure of evidence should satisfy the
principle of evidence, namely, the existence of a cut-off value that
determines evidence for and against as prescribed by the principle. Naturally,
this cut-off is 1 for the relative belief ratio. The Bayes factor is also a
valid measure of evidence and with the same cut-off. When $\Pi_{\Psi}(A)>0$
then the Bayes factor of $A$ equals $RB(A\,|\,x)/RB(A^{c}\,|\,x)$ and so can
be defined in terms of the relative belief ratio, but not conversely. Also,
$RB(A\,|\,x)>1$ iff $RB(A^{c}\,|\,x)<1$ and so the Bayes factor is not really
a comparison of the evidence for $A$ being true with the evidence for its
negation. In the continuous case, if we define the Bayes factor for $\psi$ as
a limit as in (\ref{eq1}), then this limit equals $RB_{\Psi}(\psi\,|\,x).$
Further discussion on the choice of a measure of evidence can be found in
Evans (2015) as there are other candidates beyond these two. It is important
to note, however, that the discussion of bias depends only on the principle of
evidence and is the same no matter what valid measure of evidence is used.

The following example is carried along as it illustrates a number of
things.\smallskip

\noindent\textbf{Example 1. }\textit{Location normal.}

Suppose $x=(x_{1},\ldots,x_{n})$ is i.i.d.\ $N(\mu,\sigma_{0}^{2})$ with $\pi$
a $N(\mu_{0},\tau_{0}^{2})$ prior. Then $\mu\,|\,x\sim N(\left(  n/\sigma
_{0}^{2}+1/\tau_{0}^{2}\right)  ^{-1}(n\bar{x}/\sigma_{0}^{2}+\mu_{0}/\tau
_{0}^{2}),\left(  n/\sigma_{0}^{2}+1/\tau_{0}^{2}\right)  ^{-1})$ so%
\[
RB(\mu\,|\,x)=\left(  1+\frac{n\tau_{0}^{2}}{\sigma_{0}^{2}}\right)
^{1/2}\exp\left\{
\begin{array}
[c]{c}%
-\frac{1}{2}\left(  1+\frac{\sigma_{0}^{2}}{n\tau_{0}^{2}}\right)  \left(
\frac{\sqrt{n}(\bar{x}-\mu)}{\sigma_{0}}+\frac{\sigma_{0}(\mu_{0}-\mu)}%
{\sqrt{n}\tau_{0}^{2}}\right)  ^{2}\\
+\frac{\left(  \mu-\mu_{0}\right)  ^{2}}{2\tau_{0}^{2}}%
\end{array}
\right\}  .
\]

\subsection{Bias in Hypothesis Assessment Problems}

The value $RB_{\Psi}(\psi_{\ast}\,|\,x)$ tells us if we have evidence for or
against $H_{0}:\Psi(\theta)=\psi_{\ast}$.\smallskip

\noindent\textbf{Example 1. }\textit{Location normal (continued).}

Observe that as $\tau_{0}^{2}\rightarrow\infty,$ then $RB(\mu
\,|\,x)\rightarrow\infty$ for every $\mu$ and in particular for a hypothesized
value $H_{0}=\{\mu_{\ast}\}.$ So it would appear that overwhelming evidence is
obtained for the hypothesis when the prior is very diffuse and this holds
irrespective of what the data says. Also, when the standardized value
$\sqrt{n}|\bar{x}-\mu_{\ast}|$ is fixed, then $RB(\mu_{\ast}\,|\,x)\rightarrow
\infty$ as $n\rightarrow\infty.$ This phenomenon also occurs if a Bayes factor
(which equals $RB(\mu_{\ast}\,|\,x)$ in this case) or a posterior probability
based upon a discrete prior mass at $\mu_{\ast}$ is used to assess $H_{0}.$
Accordingly all these measures lead to a sharp disagreement with the
frequentist p-value $2(1-\Phi(\sqrt{n}|\bar{x}-\mu_{\ast}|/\sigma_{0}))$ when
it is small$.$ This is the Jeffreys-Lindley paradox and it arises quite
generally.\smallskip

The Jeffreys-Lindley paradox shows that the strength of evidence cannot be
measured strictly by the size of the measure of evidence. A logical way to
assess this is to compare the evidence for $\psi_{\ast}$ with the evidence for
the other possible values for $\psi.$ The \textit{strength} of the evidence
can then be measured by%
\begin{equation}
\Pi_{\Psi}(RB_{\Psi}(\psi\,|\,x)\leq RB_{\Psi}(\psi_{\ast}\,|\,x)\,|\,x),
\label{eq2}%
\end{equation}
the posterior probability that the true value has evidence no greater than the
evidence for $\psi_{\ast}.$ So if $RB_{\Psi}(\psi_{\ast}\,|\,x)<1$ and
(\ref{eq2}) is small, then there is strong evidence against $\psi_{\ast}$
while, if $RB_{\Psi}(\psi_{\ast}\,|\,x)>1$ and (\ref{eq2}) is large, then
there is strong evidence in favor of $\psi_{\ast}.$ The inequalities
$\Pi_{\Psi}(\{\psi_{\ast}\}\,|\,x)\leq\Pi_{\Psi}(RB_{\Psi}(\psi\,|\,x)\leq
RB_{\Psi}(\psi_{\ast}\,|\,x)\,|\,x)\leq RB_{\Psi}(\psi_{\ast}\,|\,x)$ hold and
so when $RB_{\Psi}(\psi_{\ast}\,|\,x)$ is small there is strong evidence
against $\psi_{\ast}$\ and when $RB_{\Psi}(\psi_{\ast}\,|\,x)>1$ and
$\Pi_{\Psi}(\{\psi_{\ast}\}\,|\,x)$ is big, then there is strong evidence in
favor of $\psi_{\ast}.$ Note, however, that $\Pi_{\Psi}(\{\psi_{\ast
}\}\,|\,x)\approx1$ does not guarantee $RB_{\Psi}(\psi_{\ast}\,|\,x)>1$ and if
$RB_{\Psi}(\psi_{\ast}\,|\,x)<1$ this means that there is weak evidence
against $\psi_{\ast}.$ There is no reason why multiple measures of the
strength of the evidence can't be used (see the discussion in Section 2.2).
There are some issues with (\ref{eq2}) in the continuous case that require a
modification and we refer to Evans (2015) for this as the strength does not
play a key role in the discussion here. The important point is to somehow
calibrate the measure of evidence using probability to measure how strong
belief in the evidence is.\smallskip

\noindent\textbf{Example 1. }\textit{Location normal (continued).}

A simple calculation shows that, with $\sqrt{n}|\bar{x}-\mu_{\ast}|$ fixed,
then (\ref{eq2}) converges to $2(1-\Phi(\sqrt{n}|\bar{x}-\mu_{\ast}%
|/\sigma_{0}))$ as $n\tau_{0}^{2}\rightarrow\infty.$ So, if the p-value is
small, this indicates that a large value of $RB_{\Psi}(\mu_{\ast}\,|\,x)$ is
only weak evidence in favor of $\mu_{\ast}.$ It is to be noted that the
p-value $2(1-\Phi(\sqrt{n}|\bar{x}-\mu_{\ast}|/\sigma_{0}))$ is not a valid
measure of evidence as described here because there is no cut-off that
corresponds to evidence for and evidence against. So its appearance as a
measure of the strength of the evidence is not in any sense circular.

Simple algebra shows, however, that $2(1-\Phi(\sqrt{n}|\bar{x}-\mu_{\ast
}|/\sigma_{0}))-2(1-\Phi([\log(1+n\tau_{0}^{2}/\sigma_{0}^{2})+\left(
1+n\tau_{0}^{2}/\sigma_{0}^{2}\right)  ^{-1}\left(  \bar{x}-\mu_{0}\right)
^{2}/\tau_{0}^{2}]^{1/2}),$ a difference of two p-values, is a valid measure
of evidence via the cut-off 0. From this it is seen that the values of the
first p-value $2(1-\Phi(\sqrt{n}|\bar{x}-\mu_{\ast}|/\sigma_{0})$ that lead to
evidence against generally become smaller as $n\tau_{0}^{2}\rightarrow\infty.$
For example, with $n=10,\sigma_{0}^{2}=1,\mu_{\ast}=0$ and $\sqrt{n}|\bar
{x}-\mu_{\ast}|=1.96,$ then the p-value equals $0.05.$ Setting $\mu_{0}=0$ and
$\tau_{0}^{2}=1$ the second p-value equals $0.119$ and so there is evidence
against, with $\tau_{0}^{2}=10$ the the second term equals $0.032$ and with
$\tau_{0}^{2}=100$ it equals $0.009,$ so there is evidence in favor in both
cases. When $n$ increases these values become smaller as with $n=50$ the first
p-value equal to $0.05$ is always evidence in favor. Similar results are
obtained with a uniform prior on $(-m,m),$ reflecting perhaps a desire to
treat many values equivalently, as $m\rightarrow\infty$ or $n\rightarrow
\infty.$ For example, with $m=10$ and $n=10,\sigma_{0}^{2}=1,\mu_{\ast
}=0,\sqrt{n}|\bar{x}-\mu_{\ast}|=1.96,$ then the second p-value equals $0.002$
and there is evidence in favor$.$ These conclusions are similar to those found
in Berger and Selke (1987) and Berger and Delampady (1987).

It is very simple to elicit $(\mu_{0},\tau_{0}^{2})$ based on prescribing an
interval that contains the true $\mu$ with some high probability such as
$99.9\%$, taking $\mu_{0}$ to be the mid-point and so $\tau_{0}^{2}$ is
determined. There is no reason to take $\tau_{0}^{2}$ to be arbitrarily large.
But still one wonders if the choice made is inducing some kind of bias into
the problem as taking $\tau_{0}^{2}$ too large clearly does.\smallskip

Certainly default choices of priors should be avoided when possible, but even
when eliciting, how can we know if the chosen prior is inducing bias? To
assess this a numerical measure is required. The principle of evidence
suggests that \textit{bias against} $H_{0}$ is measured by
\begin{equation}
M(RB_{\Psi}(\psi_{\ast}\,|\,X)\leq1\,|\,\psi_{\ast}) \label{biasaghyp}%
\end{equation}
where $M(\cdot\,|\,\psi_{\ast})$ is the prior predictive distribution of the
data given that the hypothesis is true. So (\ref{biasaghyp}) is the prior
probability that evidence in favor of $\psi_{\ast}$ will not be obtained when
$\psi_{\ast}$ is the true value. If (\ref{biasaghyp}) is large, then there is
an \textit{a priori} bias against $H_{0}.$

For the bias in favor of $H_{0}$ it is necessary to assess if evidence against
$H_{0}$ will not be obtained with high prior probability even when $H_{0}$ is
false. One possibility is to measure \textit{bias in favor} by%
\begin{align}
&  \int_{\Psi\backslash\{\psi_{\ast}\}}M(RB_{\Psi}(\psi_{\ast}\,|\,X)\geq
1\,|\,\psi)\,\Pi_{\Psi}(d\psi)\nonumber\\
&  =M(RB_{\Psi}(\psi_{\ast}\,|\,X)\geq1)-M(RB_{\Psi}(\psi_{\ast}%
\,|\,X)\geq1\,|\,\psi_{\ast})\Pi_{\Psi}(\{\psi_{\ast}\}) \label{biasfavhyp1}%
\end{align}
which is the prior probability of not obtaining evidence against $\psi_{\ast}$
when it is false. When $\Pi_{\Psi}(\{\psi_{\ast}\})=0,$ then
(\ref{biasfavhyp1}) equals $M(RB_{\Psi}(\psi_{\ast}\,|\,X)\geq1)$ where $M$ is
the prior predictive for the data. For continuous parameters it can be argued
that it does not make sense to consider values of $\psi$ so close to
$\psi_{\ast}\ $that they are practically speaking indistinguishable. Suppose
then there is a measure of distance $d_{\Psi}$ on $\Psi$ and a value
$\delta>0$ such that, if $d_{\Psi}(\psi_{\ast},\psi)<\delta,$ then $\psi
_{\ast}$ and $\psi$ are indistinguishable in the application. The \textit{bias
in favor} of $H_{0}$ can then be measured by replacing $\Psi\backslash
\{\psi_{\ast}\}$ in (\ref{biasfavhyp1}) by $\{\psi:d_{\Psi}(\psi_{\ast}%
,\psi)\geq\delta\}$ which has upper bound bound%
\begin{equation}
\sup_{\psi:d_{\Psi}(\psi_{\ast},\psi)\geq\delta}M(RB_{\Psi}(\psi_{\ast
}\,|\,X)\geq1\,|\,\psi). \label{biasfavhyp2}%
\end{equation}
Typically $M(RB_{\Psi}(\psi_{\ast}\,|\,X)\geq1\,|\,\psi)$ decreases as $\psi$
moves away from $\psi_{\ast}$ so (\ref{biasfavhyp2}) can be computed by
finding the supremum over the set $\{\psi:d_{\Psi}(\psi_{\ast},\psi)=\delta\}$
and, when $\psi$ is real-valued and $d_{\Psi}$ is Euclidian distance, this
equals $\{\psi_{\ast}-\delta,\psi_{\ast}+\delta\}.$

It is to be noted that the measures of bias given by (\ref{biasaghyp}),
(\ref{biasfavhyp1}) and (\ref{biasfavhyp2}) do not depend on using the
relative belief ratio to measure evidence. Any valid measure of evidence will
determine the same values when the relevant cut-off is substituted for 1. It
is only (\ref{eq2}) that depends on the specific choice of the relative belief
ratio as the measure of evidence.

Under general circumstances, see Evans (2015), both biases will converge to 0
as the amount of data increases and so both biases can be controlled by
design. Clearly there is no point in reporting the results of an analysis when
there is a lot of bias unless the evidence actually contradicts the
bias.\smallskip

\noindent\textbf{Example 1. }\textit{Location normal (continued).}

Under $M(\cdot\,|\,\mu),$ then $\bar{x}\sim N(\mu,\tau_{0}^{2}+\sigma_{0}%
^{2}/n).\ $So, putting
\begin{align*}
a(\mu_{\ast},\mu_{0},\tau_{0}^{2},\sigma_{0}^{2},n)  &  =\sigma_{0}(\mu_{\ast
}-\mu_{0})/\sqrt{n}\tau_{0}^{2},\\
b(\mu_{\ast},\mu_{0},\tau_{0}^{2},\sigma_{0}^{2},n)  &  =\left\{  \left(
1+\frac{\sigma_{0}^{2}}{n\tau_{0}^{2}}\right)  \left[  \log\left(
1+\frac{n\tau_{0}^{2}}{\sigma_{0}^{2}}\right)  +\frac{\left(  \mu_{\ast}%
-\mu_{0}\right)  ^{2}}{\tau_{0}^{2}}\right]  \right\}  ^{1/2},
\end{align*}
then (\ref{biasaghyp}) is given by
\begin{align}
M(RB(\mu_{\ast}\,|\,X)\leq1\,|\,\mu_{\ast})=1-  &  \Phi\left(  a(\mu_{\ast
},\mu_{0},\tau_{0}^{2},\sigma_{0}^{2},n)+b(\mu_{\ast},\mu_{0},\tau_{0}%
^{2},\sigma_{0}^{2},n)\right)  +\nonumber\\
&  \Phi\left(  a(\mu_{\ast},\mu_{0},\tau_{0}^{2},\sigma_{0}^{2},n)-b(\mu
_{\ast},\mu_{0},\tau_{0}^{2},\sigma_{0}^{2},n)\right)  . \label{biasaglocnorm}%
\end{align}
This goes to 0 as $n\rightarrow\infty$ or as $\tau_{0}^{2}\rightarrow\infty.$
So bias against can be controlled by sample size $n$ or by the diffuseness of
the prior although, as subsequently shown, a diffuse prior induces bias in
favor. It is also the case that (\ref{biasaglocnorm}) converges to 0 when
$\mu_{0}\rightarrow\pm\infty\ $or when $\sigma_{0}/\sqrt{n}\tau_{0}$ is fixed
and $\tau_{0}\rightarrow0.$ So it would appear that using a prior with
location quite different than the hypothesized value or a prior that was much
more concentrated than the sampling distribution, can be used to lower bias
against. These are situations, however, where one can expect to have
prior-data conflict after observing the data.

\begin{figure}[pt]
\begin{center}
\includegraphics[
height=2.7882in,
width=2.7882in
]%
{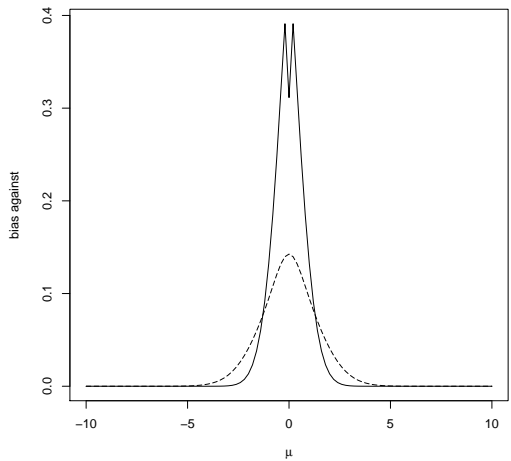}
\end{center}
\caption{Plot of bias against $H_{0}=\{\mu\}$ with a $N(0,1)$ prior (- - -)
and a $N(0,0.01)$ prior (---) with $n=5,\sigma_{0}=1.$}%
\label{fig1}%
\end{figure}

The entries in Table \ref{loctab1} record the bias against for a specific case
and illustrate that increasing $n$ does indeed reduce bias. The entries also
show that bias against can be greater when the prior is centered on the
hypothesis. Figure \ref{fig1} contains a plot of the bias against $H_{0}%
=\{\mu_{\ast}\},$ as a function of $\mu_{\ast},$ when using a $N(0,1)$ prior.
Note that the maximum bias against occurs at the mean of the prior (and equals
$0.143$) and this typically occurs when $\sigma_{0}^{2}/n\tau_{0}^{2}<1,$
namely, when the data is more concentrated than the prior. Figure \ref{fig1}
also contains a plot of the bias against when using a prior more concentrated
than the data distribution. That the bias against is maximized, as a function
of the hypothesized mean $\mu_{\ast},$ when $\mu_{\ast}$ equals the value
associated with the strongest belief under the prior seems odd. This
phenomenon arises quite often, and the mathematical explanation for this is
that the greater the amount of prior probability assigned to a value, the
harder it is for the posterior probability to increase and so it is quite
logical when considering evidence. It will be seen that this phenomenon is
very convenient for the control of bias in estimation problems and could be
used as an argument for using a prior centered on the hypothesis, although
this is not necessary as beliefs may be different.%

%TCIMACRO{\TeXButton{B}{\begin{table}[tbp] \centering}}%
%BeginExpansion
\begin{table}[tbp] \centering
%EndExpansion%
\begin{tabular}
[c]{|c|c|c|}\hline
$n$ & $\mu_{0}=1,\tau_{0}=1$ & $\mu_{0}=0,\tau_{0}=1$\\\hline
$5$ & $0.095$ & $0.143$\\
$10$ & $0.065$ & $0.104$\\
$20$ & $0.044$ & $0.074$\\
$50$ & $0.026$ & $0.045$\\
$100$ & $0.018$ & $0.031$\\\hline
\end{tabular}
\caption{Bias against for the hypothesis
$H_0=\{0\}$ with a $N(\mu_0,\tau_0^2)$ prior for different sample sizes $n$ with $\sigma_0=1$.}\label{loctab1}%
%TCIMACRO{\TeXButton{E}{\end{table}}}%
%BeginExpansion
\end{table}%
%EndExpansion

\begin{figure}
[ptb]
\begin{center}
\includegraphics[
height=2.4915in,
width=2.4915in
]%
{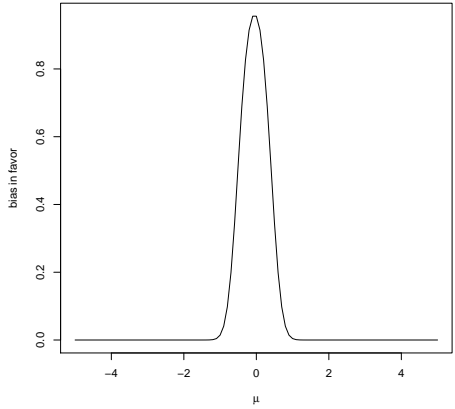}%
\caption{Plot of $M(RB(0\,|\,X)\geq1\,|\,\mu)$ when $n=20,\mu_{0}=1,\tau
_{0}=1,\sigma_{0}=1.$}%
\label{fig2}
\end{center}
\end{figure}

Now consider (\ref{biasfavhyp2}), namely, bias in favor of $H_{0}=\{\mu_{\ast
}\}.$ Putting
\[
c(\mu_{\ast},\mu,\mu_{0},\tau_{0}^{2},\sigma_{0}^{2},n)=\sqrt{n}(\mu_{\ast
}-\mu)/\sigma_{0}+a(\mu_{\ast},\mu_{0},\tau_{0}^{2},\sigma_{0}^{2},n),
\]
then (\ref{biasfavhyp2}) equals $\max M(RB(\mu_{\ast}\,|\,X)\geq
1\,|\,\mu_{\ast}\pm\delta)$ where
\begin{align}
M(RB(\mu_{\ast}\,|\,X)\geq1\,|\,\mu)=\,  &  \Phi\left(  c(\mu_{\ast},\mu
,\mu_{0},\tau_{0}^{2},\sigma_{0}^{2},n)+b(\mu_{\ast},\mu_{0},\tau_{0}%
^{2},\sigma_{0}^{2},n)\right)  -\nonumber\\
&  \Phi\left(  c(\mu_{\ast},\mu,\mu_{0},\tau_{0}^{2},\sigma_{0}^{2}%
,n)-b(\mu_{\ast},\mu_{0},\tau_{0}^{2},\sigma_{0}^{2},n)\right)
\label{biasforlocnorm}%
\end{align}
which converges to 0 as $n\rightarrow\infty$ and also as $\mu\rightarrow
\pm\infty.$ But (\ref{biasforlocnorm}) converges to 1 as $\tau_{0}%
^{2}\rightarrow\infty,$ so if the prior is too diffuse there will be bias in
favor of $\mu_{\ast}.$ So resolving the Jeffreys-Lindley paradox requires
choosing the sample size $n$, after choosing the prior, so that
(\ref{biasforlocnorm}) is suitably small. Note that choosing $\tau_{0}^{2}$
larger reduces bias against but increases bias in favor and so generally bias
cannot be avoided by choice of prior. Figure \ref{fig2} is a plot of
$M(RB(\mu_{\ast}\,|\,X)\geq1\,|\,\mu)$ for a particular case and this strictly
decreases as $\mu$ moves away from $\mu_{\ast}$.%
%TCIMACRO{\FRAME{ftbpFU}{2.4915in}{2.4915in}{0pt}{\Qcb{Plot of
%$M(RB(0\,|\,X)\geq1\,|\,\mu)$ when $n=20,\mu_{0}=1,\tau_{0}=1,\sigma_{0}=1.$}%
%}{}{fig2.eps}{\special{ language "Scientific Word";  type "GRAPHIC";
%width 2.4915in;  height 2.4915in;  depth 0pt;  original-width 6.9998in;
%original-height 6.9998in;  cropleft "0";  croptop "1";  cropright "1";
%cropbottom "0";  filename 'locnormal/fig2.eps';file-properties "XNPEU";}}}%
%BeginExpansion
%EndExpansion

In Table 2 we have recorded some specific values of the bias in favor using
(\ref{biasfavhyp1}) and using (\ref{biasfavhyp2}) where $d_{\Psi}$ is
Euclidean distance. It is seen that bias in favor can be quite serious for
small samples. When using (\ref{biasfavhyp2}) this can be mitigated by making
$\delta$ larger. For example, with $(\mu_{0},\tau_{0})=(0,1),\delta=1.0,n=20$
the bias in favor equals $0.004.$ Note, however, that $\delta$ is not chosen
to make the bias in favor small, rather it is determined in an application as
the difference from the null that is just practically important. The virtues
of determining a suitable value of $\delta$ are also readily apparent as
(\ref{biasfavhyp2}) is much smaller than (\ref{biasfavhyp1}) for larger $n.$%

%TCIMACRO{\TeXButton{B}{\begin{table}[tbp] \centering}}%
%BeginExpansion
\begin{table}[tbp] \centering
%EndExpansion%
\begin{tabular}
[c]{|c|c|c|}\hline
$n$ & $(\mu_{0},\tau_{0})=(1,1)$ & $(\mu_{0},\tau_{0})=(0,1)$\\\hline
$5$ & $0.323\,(0.871)$ & $0.451\,(0.631)$\\
$10$ & $0.259\,(0.747)$ & $0.371\,(0.516)$\\
$20$ & $0.215\,(0.519)$ & $0.299\,(0.327)$\\
$50$ & $0.153\,(0.125)$ & $0.219\,(0.062)$\\
$100$ & $0.116\,(0.006)$ & $0.168\,(0.002)$\\\hline
\end{tabular}
\caption{Bias in favor of the hypothesis
$H_0=\{0\}$ with a $N(\mu_0,\tau_0^2)$ prior  for different sample sizes $n$ with $\sigma_0=1$ using (4) (and using (5) with $\delta=0.5$).}\label{loctab2}%
%TCIMACRO{\TeXButton{E}{\end{table}}}%
%BeginExpansion
\end{table}%
%EndExpansion

A comparison of Tables \ref{loctab1} and \ref{loctab2} shows that a study
whose purpose is to demonstrate evidence in favor of $H_{0}$\ is much more
demanding than one whose purpose is to determine whether or not there is
evidence against $H_{0}.$

\subsection{Bias in Estimation Problems}

The relative belief estimate of $\psi=\Psi(\theta)$ is the value that
maximizes the measure of evidence, namely, $\psi(x)=\arg\sup RB_{\Psi}%
(\psi\,|\,x).$ It is easy to show that $RB_{\Psi}(\psi(x)\,|\,x)\geq1$ with
the inequality strict except in trivial contexts. The accuracy of this
estimate can be measured by the "size" of the \textit{plausible region
}$Pl_{\Psi}(x)=\{\psi:RB_{\Psi}(\psi\,|\,x)>1\},$ the set of values of $\psi$
that have evidence in their favor and note $\psi(x)\in Pl_{\Psi}(x).$ To say
that $\psi(x)$ is an accurate estimate, requires that $Pl_{\Psi}(x)$ be
"small", perhaps as measured by $Vol(Pl_{\Psi}(x))$ where $Vol$ is some
measure of volume, and also have high posterior content $\Pi_{\Psi}(Pl_{\Psi
}(x)\,|\,x)$ which measures the belief that the true value is in $Pl_{\Psi
}(x).$ Note that $Pl_{\Psi}(x)$ does not depend on the specific measure of
evidence chosen, in this case the relative belief ratio. Any valid estimator
must satisfy the principle of evidence and so be in $Pl_{\Psi}(x).$ It is
argued that in an estimation problem, bias is measured by various coverage
probabilities for the plausible region.

Note too that if there is evidence in favor of $H_{0}:\Psi(\theta)=\psi_{\ast
},$ then $\psi_{\ast}\in Pl_{\Psi}(x)$ and so represents the natural estimate
of $\psi$ provided there was a clear reason for assessing the evidence for
this value. The strength of the evidence in favor of $\psi_{\ast}$ can then
also be measured by the size of $Pl_{\Psi}(x).$ Similarly, if evidence against
$H_{0}$ is obtained then $\psi_{\ast}\in Im_{\Psi}=\{\psi:RB_{\Psi}%
(\psi\,|\,x)<1\},$ the \textit{implausible region }and then there is strong
evidence against $H_{0}$ provided $Pl_{\Psi}(x)$ has small volume and large
posterior probability. A virtue of this approach to measuring the strength of
the evidence is that it does not depend upon using the relative belief ratio
to measure evidence.

\subsubsection{Bias Against}

The prior probability that the plausible region does not cover the true value
measures bias against when estimating $\psi.$ For if this probability is
large, then the estimate and the plausible region are \textit{a priori} likely
to be misleading as to the true value. The prior probability that $Pl_{\Psi
}(x)$ doesn't contain $\psi=\Psi(\theta)$ when $\theta\sim\Pi,X\sim P_{\theta
}$ is
\begin{equation}
E_{\Pi_{\Psi}}\left(  M(\psi\notin Pl_{\Psi}(X)\,|\,\psi)\right)
=E_{\Pi_{\Psi}}(M(RB_{\Psi}(\psi\,|\,X)\leq1\,|\,\psi)) \label{biasagest1}%
\end{equation}
which is also the average bias against over all hypothesis testing problems
$H_{0}:\Psi(\theta)=\psi.$ Note that $1-E_{\Pi_{\Psi}}\left(  M(\psi\notin
Pl_{\Psi}(X)\,|\,\psi)\right)  =E_{\Pi_{\Psi}}\left(  M(\psi\in Pl_{\Psi
}(X)\,|\,\psi)\right)  $\newline$=E_{M}\left(  \Pi_{\Psi}(Pl_{\Psi
}(X)\,|\,X)\right)  $ which is the prior coverage probability of $Pl_{\Psi}$.
Also,
\begin{equation}
\sup_{\psi}M(\psi\notin Pl_{\Psi}(X)\,|\,\psi)=\sup_{\psi}M(RB_{\Psi}%
(\psi\,|\,X)\leq1\,|\,\psi), \label{biasagest2}%
\end{equation}
is an upper bound on (\ref{biasagest1}). Therefore, controlling
(\ref{biasagest2}) controls the bias against in estimation and all hypothesis
assessment problems involving $\psi$. Also $1-\sup_{\psi}M(\psi\notin
Pl_{\Psi}(X)\,|\,\psi)=\inf_{\psi}M(\psi\in Pl_{\Psi}(X)\,|\,\psi)\leq
E_{M}\left(  \Pi_{\Psi}(Pl_{\Psi}(X)\,|\,X)\right)  $ so using
(\ref{biasagest2}) implies lower bounds for the coverage probability and for
the expected posterior content of the plausible region. In general, both
(\ref{biasagest1}) and (\ref{biasagest2}) converge to 0 with increasing
amounts of data. So it is possible to control for bias against in estimation
problems by design.\smallskip

\noindent\textbf{Example 1. }\textit{Location normal (continued).}

The value of $M(RB(\mu\,|\,X)\leq1\,|\,\mu)$ is given in (\ref{biasaglocnorm})
and examples are plotted in Figure \ref{fig1}. When $\mu\sim N(\mu_{0}%
,\tau_{0}^{2})$ then $z=(\mu-\mu_{0})/\tau_{0}\sim N(0,1)$ so
\begin{align*}
&  E_{\Pi}\left(  M(RB(\mu\,|\,X)\leq1\,|\,\mu)\right) \\
&  =1-E\left[
\begin{array}
[c]{c}%
\Phi\left(  \frac{\sigma_{0}}{\sqrt{n}\tau_{0}}Z+\left\{  \left(
1+\frac{\sigma_{0}^{2}}{n\tau_{0}^{2}}\right)  \left[  \log\left(
1+\frac{n\tau_{0}^{2}}{\sigma_{0}^{2}}\right)  +Z^{2}\right]  \right\}
^{1/2}\right)  +\\
\Phi\left(  \frac{\sigma_{0}}{\sqrt{n}\tau_{0}}Z-\left\{  \left(
1+\frac{\sigma_{0}^{2}}{n\tau_{0}^{2}}\right)  \left[  \log\left(
1+\frac{n\tau_{0}^{2}}{\sigma_{0}^{2}}\right)  +Z^{2}\right]  \right\}
^{1/2}\right)
\end{array}
\right]
\end{align*}
which is notably independent of the prior mean $\mu_{0}$. The dominated
convergence theorem implies $E_{\Pi}\left(  M(RB(\mu\,|\,X)\leq1\,|\,\mu
)\right)  \rightarrow0$ as $n\rightarrow\infty$ or as $\tau_{0}^{2}%
\rightarrow\infty.$ So provided $n\tau_{0}^{2}/\sigma_{0}^{2}$ is large
enough, there is no estimation bias against. Table \ref{loctab3} illustrates
some values of this bias measure. Subtracting the probabilities in Table
\ref{loctab3} from 1 gives the prior probability that the plausible region
covers the true value and the expected posterior content of the plausible
region. So when $n=20,\tau_{0}=1,$ the prior probability of the plausible
region containing the true value is $1-0.051=0.949$ so $Pl(x)$ is a $0.949$
Bayesian confidence interval for $\mu.$%

%TCIMACRO{\TeXButton{B}{\begin{table}[tbp] \centering}}%
%BeginExpansion
\begin{table}[tbp] \centering
%EndExpansion%
\begin{tabular}
[c]{|c|c|c|}\hline
$n$ & $\tau_{0}=1$ & $\tau_{0}=0.5$\\\hline
$5$ & $0.107$ & $0.193$\\
$10$ & $0.075$ & $0.146$\\
$20$ & $0.051$ & $0.107$\\
$50$ & $0.031$ & $0.067$\\
$100$ & $0.021$ & $0.046$\\\hline
\end{tabular}
\caption{Average bias against $H_0={0}$  when using a$N(0,\tau_0^2)$ prior for different sample sizes
$n$.}\label{loctab3}%
%TCIMACRO{\TeXButton{E}{\end{table}}}%
%BeginExpansion
\end{table}%
%EndExpansion

To use (\ref{biasagest2}) it is necessary to maximize $M(RB(\mu\,|\,X)\leq
1\,|\,\mu)$ as a function of $\mu$ and it is seen that, at least when the
prior is not overly concentrated, that this maximum occurs at $\mu_{0}.$
Figure \ref{fig1} shows that when using the $N(0,1)$ prior the maximum occurs
at $\mu=0$ when $n=5$ and from the second column of Table \ref{loctab1}, the
maximum equals $0.143$. The average bias against is given by $0.107,$ as
recorded in Table \ref{loctab3}. Note that the maximum also occurs at $\mu=0$
for the other values of $n$ recorded in Table \ref{loctab1}.

\subsubsection{Bias in\ Favor}

\begin{figure}
[ptb]
\begin{center}
\includegraphics[
height=2.783in,
width=2.783in
]%
{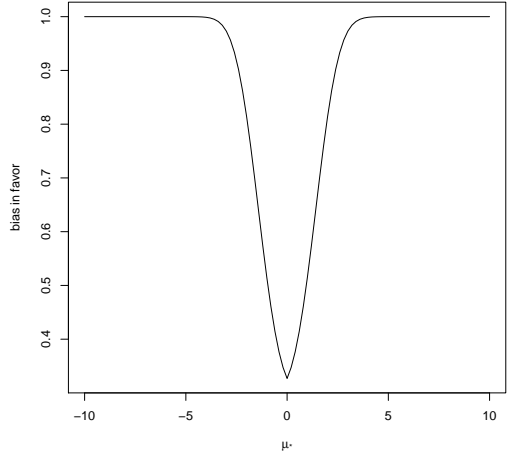}%
\caption{Bias in favor of\ $\mu$\ maximized over $\mu\pm\delta$ based on a
$N(0,1)$ prior and $\sigma_{0}=1,n=20,\delta=0.5.$}%
\label{fig3}
\end{center}
\end{figure}

Bias in favor occurs when the prior probability that $Im_{\Psi}$ does not
cover a false value is large, namely, when%
\begin{align}
&  \int_{\Psi}\int_{\Psi\backslash\{\psi_{\ast}\}}M(\psi_{\ast}\notin
Im_{\Psi}(X)\,|\,\psi)\,\Pi_{\Psi}(d\psi)\,\Pi_{\Psi}(d\psi_{\ast})\nonumber\\
&  =\int_{\Psi}\int_{\Psi\backslash\{\psi_{\ast}\}}M(RB_{\Psi}(\psi_{\ast
}\,|\,X)\geq1\,|\,\psi)\,\Pi_{\Psi}(d\psi)\,\Pi_{\Psi}(d\psi_{\ast})
\label{biasfavest1}%
\end{align}
is large as this would seem to imply that the plausible region will cover a
randomly selected false value from the prior with high prior probability. Note
that (\ref{biasfavest1}) is the prior mean of (\ref{biasfavhyp1}) and in the
continuous case equals $\int_{\Psi}M(\psi_{\ast}\notin Im_{\Psi}%
(X))\,\Pi_{\Psi}(d\psi_{\ast})$. As previously discussed, however, it often
doesn't make sense to distinguish values of $\psi$ that are close to
$\psi_{\ast}.$ The bias in favor for estimation can then be measured by%
\begin{align}
&  E_{\Pi_{\Psi}}\left(  \sup_{\psi:d_{\Psi}(\psi,\psi_{\ast})\geq\delta
}M(\psi_{\ast}\notin Im_{\Psi}(X)\,|\,\psi)\right) \nonumber\\
&  =E_{\Pi_{\Psi}}\left(  \sup_{\psi:d_{\Psi}(\psi,\psi_{\ast})\geq\delta
}M(RB_{\Psi}(\psi_{\ast}\,|\,X)\geq1\,|\,\psi)\right)  . \label{biasfavest2}%
\end{align}
An upper bound on (\ref{biasfavest2}) is commonly equal to 1 as illustrated in
Figure \ref{fig3} and so is not useful.

It is the size and posterior content of $Pl_{\Psi}(x)$ that provides a measure
of the accuracy of the estimate $\psi(x).$ As discussed in Section 2.2.1 the
\textit{a priori} expected posterior content of $Pl_{\Psi}(x)$ can be
controlled by bias against. The\textit{ a priori} expected volume of
$Pl_{\Psi}(x)$ satisfies
\begin{equation}
E_{M}\left(  Vol(Pl_{\Psi}(X))\right)  =\int_{\Psi}\int_{\Psi}M(\psi_{\ast}\in
Pl_{\Psi}(X)\,|\,\psi)\,\Pi_{\Psi}(d\psi)\,Vol(d\psi_{\ast}). \label{eq8}%
\end{equation}
Notice that when $\Pi_{\Psi}(\{\psi\})=0$ for every $\psi,$ this can be
interpreted as a kind of average of the prior probabilities of the plausible
region covering a false value.$\smallskip$

\noindent\textbf{Example 1. }\textit{Location normal (continued).}

It follows from (\ref{biasforlocnorm}) that
\begin{align*}
&  \sup M(RB(\mu_{\ast}\,|\,X)\geq1\,|\,\mu_{\ast}\pm\delta)\\
=  &  \sup\left\{
\begin{array}
[c]{c}%
\Phi\left(  c(\mu_{\ast},\mu_{\ast}\pm\delta,\mu_{0},\tau_{0}^{2},\sigma
_{0}^{2},n)+b(\mu_{\ast},\mu_{0},\tau_{0}^{2},\sigma_{0}^{2},n)\right)  -\\
\Phi\left(  c(\mu_{\ast},\mu_{\ast}\pm\delta,\mu_{0},\tau_{0}^{2},\sigma
_{0}^{2},n)-b(\mu_{\ast},\mu_{0},\tau_{0}^{2},\sigma_{0}^{2},n)\right)
\end{array}
\right\}
\end{align*}
Note that as $\mu_{\ast}\rightarrow\pm\infty$, then $M(RB(\mu_{\ast
}\,|\,X)\geq1\,|\,\mu_{\ast}\pm\delta)\rightarrow1$ when $n\tau_{0}^{2}%
/\sigma_{0}^{2}>1,$ see Figure \ref{fig3}, and converges to 0 if $n\tau
_{0}^{2}/\sigma_{0}^{2}<1,$ so it would appear that the better circumstance
for guarding against bias in favor is when the prior is putting in more
information than the data. As previously noted, however, this is a situation
where we might expect prior data-conflict to arise and, except in exceptional
circumstances should be avoided.%
%TCIMACRO{\FRAME{ftbpFU}{2.783in}{2.783in}{0pt}{\Qcb{Bias in favor of\ $\mu
%$\ maximized over $\mu\pm\delta$ based on a $N(0,1)$ prior and $\sigma
%_{0}=1,n=20,\delta=0.5.$}}{}{fig3.eps}{\special{ language "Scientific Word";
%type "GRAPHIC";  maintain-aspect-ratio TRUE;  display "USEDEF";
%valid_file "F";  width 2.783in;  height 2.783in;  depth 0pt;
%original-width 6.9998in;  original-height 6.9998in;  cropleft "0";
%filename 'locnormal/fig3.eps';file-properties "XNPEU";}}}%
%BeginExpansion
%EndExpansion
Table \ref{loctab5} contains values of (\ref{biasforlocnorm}) for this
situation with different values of $\delta$.%

%TCIMACRO{\TeXButton{B}{\begin{table}[tbp] \centering}}%
%BeginExpansion
\begin{table}[tbp] \centering
%EndExpansion%
\begin{tabular}
[c]{|c|c|c|}\hline
$n$ & $(\mu_{0},\tau_{0})=(0,1),\delta=1.0$ & $(\mu_{0},\tau_{0}%
)=(0,1),\delta=0.5$\\\hline
$5$ & $0.451$ & $0.798$\\
$10$ & $0.185$ & $0.690$\\
$20$ & $0.025$ & $0.486$\\
$50$ & $0.000$ & $0.131$\\
$100$ & $0.000$ & $0.009$\\\hline
\end{tabular}
\caption{Average bias in favor for estimation based on (11) when using a $N(0,\tau_0^2)$ prior for different sample sizes
$n$ and difference $\delta$. }\label{loctab5}%
%TCIMACRO{\TeXButton{E}{\end{table}}}%
%BeginExpansion
\end{table}%
%EndExpansion

Some elementary calculations give $Pl(x)=\bar{x}\pm w(\bar{x},n,\sigma_{0}%
^{2},\mu_{0},\tau_{0}^{2})$ with
\begin{align*}
&  w(\bar{x},n,\sigma_{0}^{2},\mu_{0},\tau_{0}^{2})\\
&  =\frac{\sigma_{0}}{\sqrt{n}}\left(  1+\frac{n\tau_{0}^{2}}{\sigma_{0}^{2}%
}\right)  ^{-1/2}\left\{  \left(  1+\frac{n\tau_{0}^{2}}{\sigma_{0}^{2}%
}\right)  \log\left(  1+\frac{n\tau_{0}^{2}}{\sigma_{0}^{2}}\right)  +\left(
\frac{\bar{x}-\mu_{0}}{\sigma_{0}/\sqrt{n}}\right)  ^{2}\right\}  ^{1/2}%
\end{align*}
where $z=\sqrt{n}(\bar{x}-\mu_{0})/\sigma_{0}\sim N(0,1$) under $M.$ It is
notable that the prior distribution of the width is independent of the prior
mean. Table \ref{loctab3.1} contains some expected half-widths together with
the coverage probabilities of $Pl(x).$%

%TCIMACRO{\TeXButton{B}{\begin{table}[tbp] \centering}}%
%BeginExpansion
\begin{table}[tbp] \centering
%EndExpansion%
\begin{tabular}
[c]{|c|c|c|}\hline
$n$ & $\tau_{0}=1$ & $\tau_{0}=0.5$\\\hline
$5$ & $0.625\,(0.893)$ & $0.491\,(0.807)$\\
$10$ & $0.499\,(0.925)$ & $0.389\,(0.854)$\\
$20$ & $0.393\,(0.949)$ & $0.312\,(0.893)$\\
$50$ & $0.281\,(0.969)$ & $0.231\,(0.933)$\\
$100$ & $0.215\,(0.979\,)$ & $0.181\,(0.954)$\\\hline
\end{tabular}
\caption{Expected  half-widths (coverages) of the plausible interval when using a$N(\mu_0,\tau_0^2)$ prior for different sample sizes
$n$.}\label{loctab3.1}%
%TCIMACRO{\TeXButton{E}{\end{table}}}%
%BeginExpansion
\end{table}%
%EndExpansion

\section{Frequentist and Optimal Properties}

Consider now the bias against $H_{0}=\{\psi_{\ast}\},$ namely, $M(RB_{\Psi
}(\psi_{\ast}\,|\,X)\leq1\,|\,\psi_{\ast}).$ If we repeatedly generate
$\theta\sim\pi(\cdot\,|\,\psi_{\ast}),X\sim f_{\theta},$ then this probability
is the long-run proportion of times that $RB_{\Psi}(\psi_{\ast}\,|\,X)\leq1.$
This frequentist interpretation depends on the conditional prior $\pi
(\cdot\,|\,\psi_{\ast})\ $and when $\Psi(\theta)=\theta,$ so there are no
nuisance parameters, this is a "pure" frequentist probability. Even in the
latter case there is some dependence on the prior, however, as $RB(\theta
_{\ast}\,|\,x)=f_{\theta_{\ast}}(x)/m(x)$ so $x$ satisfies $RB_{\Psi}%
(\theta_{\ast}\,|\,x)\leq1$ iff $f_{\theta_{\ast}}(x)\leq m(x)$ where
$m(x)=\int_{\Theta}f_{\theta}(x)\,\Pi(d\theta).$ So in general the region
$\{x:RB_{\Psi}(\psi_{\ast}\,|\,x)\leq1\}$ depends on $\pi$ but the probability
$M(RB_{\Psi}(\psi_{\ast}\,|\,X)\leq1\,|\,\psi_{\ast})$ depends only on the
conditional prior predictive given $\Psi(\theta)=\psi_{\ast},$ namely,
$m(x\,|\,\psi_{\ast})=\int_{\Theta}f_{\theta}(x)\,\Pi(d\theta\,|\,\psi_{\ast
}),$ and not on the marginal prior $\pi_{\Psi}$ on $\psi.$ We refer to
probabilities that depend only on $M(\cdot\,\,|\,\psi_{\ast})$ as frequentist,
for example, coverage probabilities are called confidences, and those that
depend on the full prior $\pi$ as Bayesian confidences. The frequentist label
is similar to use of the confidence terminology when dealing with random
effects models as nuisance parameters have been integrated out.

Suppose now that some other general rule, not necessarily the principle of
evidence, is used to determine whether there is evidence for or against
$\psi_{\ast}$\ and this leads to the set $D(\psi_{\ast})\subset\mathcal{X}$ as
those data sets that do not give evidence in favor of $H_{0}=\{\psi_{\ast}\}.$
The rules of potential interest will satisfy $M(D(\psi_{\ast})\,|\,\psi_{\ast
})\leq M(RB_{\Psi}(\psi_{\ast}\,|\,X)\leq1\,|\,\psi_{\ast})$ since this
implies better performance \textit{a priori} in terms of identifying when data
has evidence in favor of $H_{0}$ via the set $D^{c}(\psi_{\ast})$ than the
principal of evidence$.$ For example, $D(\psi_{\ast})=\{x:RB_{\Psi}(\psi
_{\ast}\,|\,x)\leq q\}$ for some $q<1$ satisfies this but note that a value
satisfying $q<RB_{\Psi}(\psi_{\ast}\,|\,x)\leq1$ violates the principle of
evidence if it is claimed there is evidence in favor of $\psi_{\ast}$. Putting
$R(\psi_{\ast})=\{x:RB_{\Psi}(\psi_{\ast}\,|\,x)\leq1\}$ leads to the
following result.\smallskip

\noindent\textbf{Theorem 1.} (i) The prior probability $M(D(\psi_{\ast}))$ is
maximized among all $D(\psi_{\ast})\subset\mathcal{X}$ satisfying
$M(D(\psi_{\ast})\,|\,\psi_{\ast})\leq M(R(\psi_{\ast})\,|\,\psi_{\ast})$ by
$D(\psi_{\ast})=R(\psi_{\ast}).$ (ii) If $\Pi_{\Psi}(\{\psi_{\ast}\})=0,$ then
$R(\psi_{\ast})$ maximizes the prior probability of not obtaining evidence in
favor of $\psi_{\ast}$ when it is false and otherwise maximizes this
probability among all rules satisfying $M(D(\psi_{\ast})\,|\,\psi_{\ast
})=M(R(\psi_{\ast})\,|\,\psi_{\ast}).$\smallskip

\noindent When $\Pi_{\Psi}(\{\psi_{\ast}\})\neq0,$ rules may exist having
greater prior probability of not getting evidence in favor of $\psi_{\ast}$
when it is false but the price paid for this is the violation of the principle
of evidence. Also, when comparing rules based on their ability to distinguish
falsity it only seems fair that the rules perform the same under the truth. So
Theorem 1 is a general optimality result for the principle of evidence applied
to hypothesis assessment when considering bias against.

Now consider $C(x)=\{\psi:x\notin D(\psi)\}$ which is the set of $\psi$ values
for which there is evidence in their favor after observing $x$ according to
some alternative evidence rule. Since $M(\psi_{\ast}\notin C(X)\,|\,\psi
)=M(D(\psi_{\ast}))\,|\,\psi),$ then $E_{\Pi_{\Psi}}\left(  M(\psi\in
C(X)\,)\,|\,\psi)\right)  =1-E_{\Pi_{\Psi}}\left(  M(\psi\notin C(X)\,|\,\psi
)\right)  =1-E_{\Pi_{\Psi}}\left(  M(D(\psi)\,|\,\psi)\right)  $\newline%
$\geq1-E_{\Pi_{\Psi}}\left(  M(R(\psi)\,|\,\psi)\right)  =E_{\Pi_{\Psi}%
}\left(  M(\psi\in Pl_{\Psi}(X)\,)\,|\,\psi)\right)  $ and so the Bayesian
coverage of $C$ is at least as large as that of $Pl_{\Psi}$ and so represents
a viable alternative to using $Pl_{\Psi}.$

The following establishes an optimality result for $Pl_{\Psi}$.\smallskip

\noindent\textbf{Theorem 2.} (i) The prior probability that the region $C$
doesn't cover a value $\psi_{\ast}$ generated from the prior, namely,
$E_{\Pi_{\Psi}}(M(\psi_{\ast}\notin C(X))),$ is maximized among all regions
satisfying $M(\psi_{\ast}\notin C(X)\,|\,\psi_{\ast})\leq M(\psi_{\ast}\notin
Pl_{\Psi}(X)\,|\,\psi_{\ast})$ for every $\,\psi_{\ast},$ by $C=Pl_{\Psi}.$
(ii) If $\Pi_{\Psi}(\{\psi_{\ast}\})=0$ for all $\psi_{\ast},$ then $Pl_{\Psi
}$ maximizes the prior probability of not covering a false value and otherwise
maximizes this probability among all $C$ satisfying $M(\psi_{\ast}\notin
C(X)\,|\,\psi_{\ast})=M(\psi_{\ast}\notin Pl_{\Psi}(X)\,|\,\psi_{\ast})$ for
all $\psi_{\ast}.$\smallskip

\noindent Again when $\Pi_{\Psi}(\{\psi_{\ast}\})\neq0$ the existence of a
region with better properties with respect to not covering false values than
$Pl_{\Psi}$ can't be ruled out but, when considering such a property, it seems
only fair to compare regions with the same coverage probability and in that
case $Pl_{\Psi}$ is optimal. So Theorem 2 is also a general optimality result
for the principle of evidence applied to estimation when considering bias
against. Also, if there is a value $\psi_{0}=\arg\inf_{\psi}M(\psi\in
Pl_{\Psi}(X)\,)\,|\,\psi),$ then $\gamma_{0}=M(\psi_{0}\in Pl_{\Psi
}(X)\,)\,|\,\psi_{0})$ serves as a lower bound on the coverage probabilities,
and thus $Pl_{\Psi}$ is a $\gamma_{0}$-confidence region for $\psi$ and this
is a pure frequentist $\gamma_{0}$-confidence region when $\Psi(\theta
)=\theta.$ Since $M(\psi\in Pl_{\Psi}(X)\,)\,|\,\psi)=1-M(\psi\notin Pl_{\Psi
}(X)\,)\,|\,\psi)=1-M(R(\psi_{\ast})\,|\,\psi),$ then Example 1 shows that it
is reasonable to expect that such a $\psi_{0}$ exists.

The principle of evidence leads to the following satisfying properties which
connect the concept of bias as discussed here with the frequentist
concept..\smallskip

\noindent\textbf{Theorem 3.} (i) Using the principle of evidence, the prior
probability of getting evidence in favor of $\psi_{\ast}$ when it is true is
greater than or equal to the prior probability of getting evidence in favor of
$\psi_{\ast}$ given that $\psi_{\ast}$ is false. (ii) The prior probability of
$Pl_{\Psi}$\ covering the true value is always greater than or equal to the
prior probability of $Pl_{\Psi}$ covering a false value.\smallskip

\noindent The properties stated in Theorem 3 are similar to a property called
unbiasedness for frequentist procedures.\ For example, a test is unbiased if
the probability of rejecting a null is always larger when it is false than
when it is true and a confidence region is unbiased if the probability of
covering the true value is always greater than the probability of covering a
false value. While the inferences discussed here are "unbiased" in this
generalized sense, they could still be biased against or in favor in the
practical sense of this paper, as it is the amount of data that controls this.

Now consider bias in favor and suppose there is an alternative
characterization of evidence that leads to the region $E(\psi_{\ast})$
consisting of all data sets that do not lead to evidence against $\psi_{\ast
}.$ Putting $A(\psi_{\ast})=\{x:RB_{\Psi}(\psi_{\ast}\,|\,x)\geq1,$ we
restrict attention to regions satisfying $M(E(\psi_{\ast})\,|\,\psi_{\ast
})\geq M(A(\psi_{\ast})\,|\,\psi_{\ast}).$ Using (\ref{biasfavhyp1}) to
measure bias in leads to the following results.\smallskip

\noindent\textbf{Theorem 4.} (i) The prior probability $M(E(\psi_{\ast}))$ is
minimized among all $E(\psi_{\ast})\subset\mathcal{X}$ satisfying
$M(E(\psi_{\ast})\,|\,\psi_{\ast})\geq M(A(\psi_{\ast})\,|\,\psi_{\ast})$ by
$E(\psi_{\ast})=A(\psi_{\ast}).$ (ii) If $\Pi_{\Psi}(\{\psi_{\ast}\})=0,$ then
the set $A(\psi_{\ast})$ minimizes the prior probability of not obtaining
evidence against $\psi_{\ast}$ when it is false and otherwise minimizes this
probability among all rules satisfying $M(E(\psi_{\ast})\,|\,\psi_{\ast
})=M(A(\psi_{\ast})\,|\,\psi_{\ast}).$\smallskip

\noindent\textbf{Theorem 5.} (i) The prior probability region $C$ covers a
value $\psi_{\ast}$ generated from the prior, namely, $E_{\Pi_{\Psi}}%
(M(\psi_{\ast}\in C(X))),$ is minimized among all regions satisfying
$M(\psi_{\ast}\in C(X)\,|\,\psi_{\ast})\geq M(\psi_{\ast}\in Pl_{\Psi
}(X)\,|\,\psi_{\ast})$ for every $\,\psi_{\ast},$ by $C=Pl_{\Psi}.$ (ii) If
$\Pi_{\Psi}(\{\psi_{\ast}\})=0$ for all $\psi_{\ast},$ then $Pl_{\Psi}$
minimizes the prior probability of covering a false value and otherwise
minimizes this probability among all rules satisfying $M(\psi_{\ast}\in
C(X)\,|\,\psi_{\ast})=M(\psi_{\ast}\in Pl_{\Psi}(X)\,|\,\psi_{\ast})$ for all
$\psi_{\ast}.$\smallskip

\noindent So Theorems 4 and 5 are optimality results for the principle of
evidence when considering bias in favor.

Clearly the bias against $H_{0}$ is playing a role similar to size in
frequentist statistics and the bias in favor is playing a role similar to
power. A study that found evidence against $H_{0},$ but had a high bias
against, or a study that found evidence in favor of $H_{0}$ but had a high
bias in favor, could not be considered to be of high quality. Similarly, a
study concerned with estimating a quantity of interest could not be considered
of high quality if there is high bias against or in favor. There are some
circumstances, however, where some bias is perhaps not an issue. For example,
in a situation where sparsity is to be expected, then allowing for high bias
in favor of certain hypotheses accompanied by low bias against, may be
tolerable although this does reduce the reliability of any hypotheses where
evidence is found in favor.

\section{Examples}

A number of examples are now considered.\smallskip

\noindent\textbf{Example 2.} \textit{Binomial.}

Suppose $x=(x_{1},\ldots,x_{n})$ is a sample from the Bernoulli$(\theta)$ with
$\theta\in\lbrack0,1]$ unknown so $n\bar{x}\sim$ binomial$(n,\theta)$ and
interest is in $\theta.$ For the prior let $\theta\sim$ beta$(\alpha_{0}%
,\beta_{0})$ where the hyperparameters are elicited as in, for example, Evans,
Guttman and Li (2017), so $\theta\,|\,n\bar{x}\sim$ beta$(\alpha_{0}+n\bar
{x},\beta_{0}+n(1-\bar{x})).$ Then
\[
RB(\theta\,|\,n\bar{x})=\frac{\Gamma(\alpha_{0}+\beta_{0}+n)}{\Gamma
(\alpha_{0}+n\bar{x})\Gamma(\beta_{0}+n(1-\bar{x}))}\frac{\Gamma(\alpha
_{0})\Gamma(\beta_{0})}{\Gamma(\alpha_{0}+\beta_{0})}\theta^{n\bar{x}%
}(1-\theta)^{n(1-\bar{x})}%
\]
is unimodal with mode at $\bar{x},$ so $Pl(x)$ is an interval containing
$\bar{x}.$ Note that $M(\cdot\,|\,\theta)$ is the binomial$(n,\theta)$
probability measure and the bias against $\theta$ is given by $M(RB(\theta
\,|\,n\bar{x})\leq1\,|\,\theta)$ while the bias in favor of $\theta$, using
(\ref{biasfavhyp2}), is given by $\max M(RB(\theta\,|\,n\bar{x})\geq
1\,|\,\theta\pm\delta)$for $\theta\in\lbrack\delta,1-\delta].$

\begin{figure}
[p]
\begin{center}
\includegraphics[
height=2.5071in,
width=2.5071in
]%
{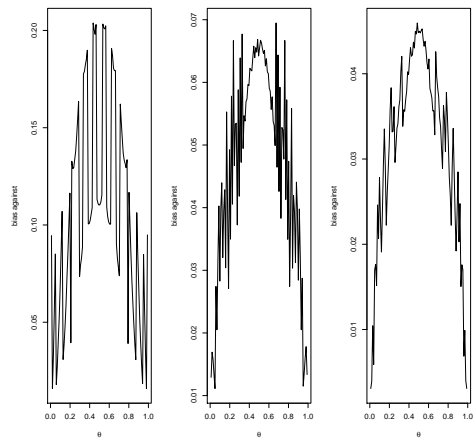}%
\caption{Plots of bias against as a function of $\theta$ for $n=10,50$ and
$100$ when using beta$(1,1)$ prior.}%
\label{binomfig1}%
\end{center}
\end{figure}

\begin{figure}
[p]
\begin{center}
\includegraphics[
height=2.5918in,
width=2.5918in
]%
{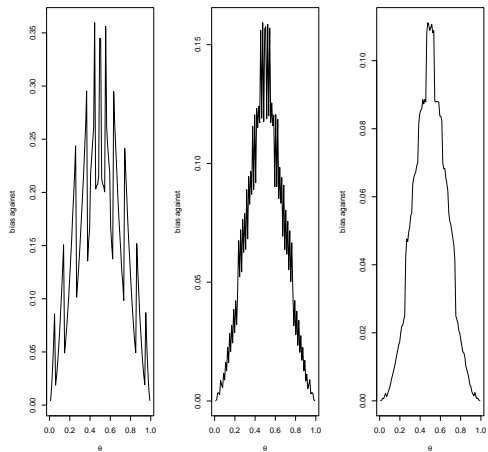}%
\caption{Plots of bias against as a function of $\theta$ for $n=10,50$ and
$100$ when using beta$(5,5)$ prior.}%
\label{binomfig2}%
\end{center}
\end{figure}

Consider first the prior given by $(\alpha_{0},\beta_{0})=(1,1).$ Figure
\ref{binomfig1} gives the plots of the bias against for $n=10$ (max. = $0.21$,
average = $0.11$), $n=50$ (max.= $0.07$, average = $0.05$) and $n=100$ (max. =
$0.05$, average = $0.03$). Therefore, when $n=10,$ then $Pl(x)$ is a
$0.79$-confidence interval for $\theta,$ when $n=50$ it is a $0.93$-confidence
interval for $\theta$ and when $n=100$ it is a $0.95$-confidence interval for
$\theta.$ For the informative prior given by $(\alpha_{0},\beta_{0})=(5,5)$,
Figure \ref{binomfig2} gives the plots of the bias against for $n=10$ (max. =
$0.36$, average = $0.21$), $n=50$ (max. = $0.16$, average = $0.10$) and
$n=100$ (max. = $0.11$, average = $0.07$). So when $n=10$ then $Pl(x)$ is a
$0.64$-confidence interval for $\theta,$ when $n=50$ it is a $0.84$-confidence
interval for $\theta$ and when $n=100$ it is a $0.93$-confidence interval for
$\theta.$ One feature immediately stands out, namely, when using a more
informative prior the bias against increases. As previously explained this
phenomenon occurs because when the prior probability of $\theta$ is small, it
is much easier to obtain evidence in favor than when the prior probability of
$\theta$ is large.%
%TCIMACRO{\FRAME{fpFU}{2.5071in}{2.5071in}{0pt}{\Qcb{Plots of bias against as a
%function of $\theta$ for $n=10,50$ and $100$ when using beta$(1,1)$ prior.}%
%}{\Qlb{binomfig1}}{binomfig1.eps}{\special{ language "Scientific Word";
%type "GRAPHIC";  maintain-aspect-ratio TRUE;  display "USEDEF";
%valid_file "F";  width 2.5071in;  height 2.5071in;  depth 0pt;
%original-width 6.9998in;  original-height 6.9998in;  cropleft "0";
%filename 'binomialexample/binomfig1.eps';file-properties "XNPEU";}}}%
%BeginExpansion
%EndExpansion%
%TCIMACRO{\FRAME{fpFU}{2.5918in}{2.5918in}{0pt}{\Qcb{Plots of bias against as a
%function of $\theta$ for $n=10,50$ and $100$ when using beta$(5,5)$ prior.}%
%}{\Qlb{binomfig2}}{binomfig2.eps}{\special{ language "Scientific Word";
%type "GRAPHIC";  maintain-aspect-ratio TRUE;  display "USEDEF";
%valid_file "F";  width 2.5918in;  height 2.5918in;  depth 0pt;
%original-width 6.9998in;  original-height 6.9998in;  cropleft "0";
%filename 'binomialexample/binomfig2.eps';file-properties "XNPEU";}}}%
%BeginExpansion
%EndExpansion

\begin{figure}
[p]
\begin{center}
\includegraphics[
height=2.2035in,
width=2.2035in
]%
{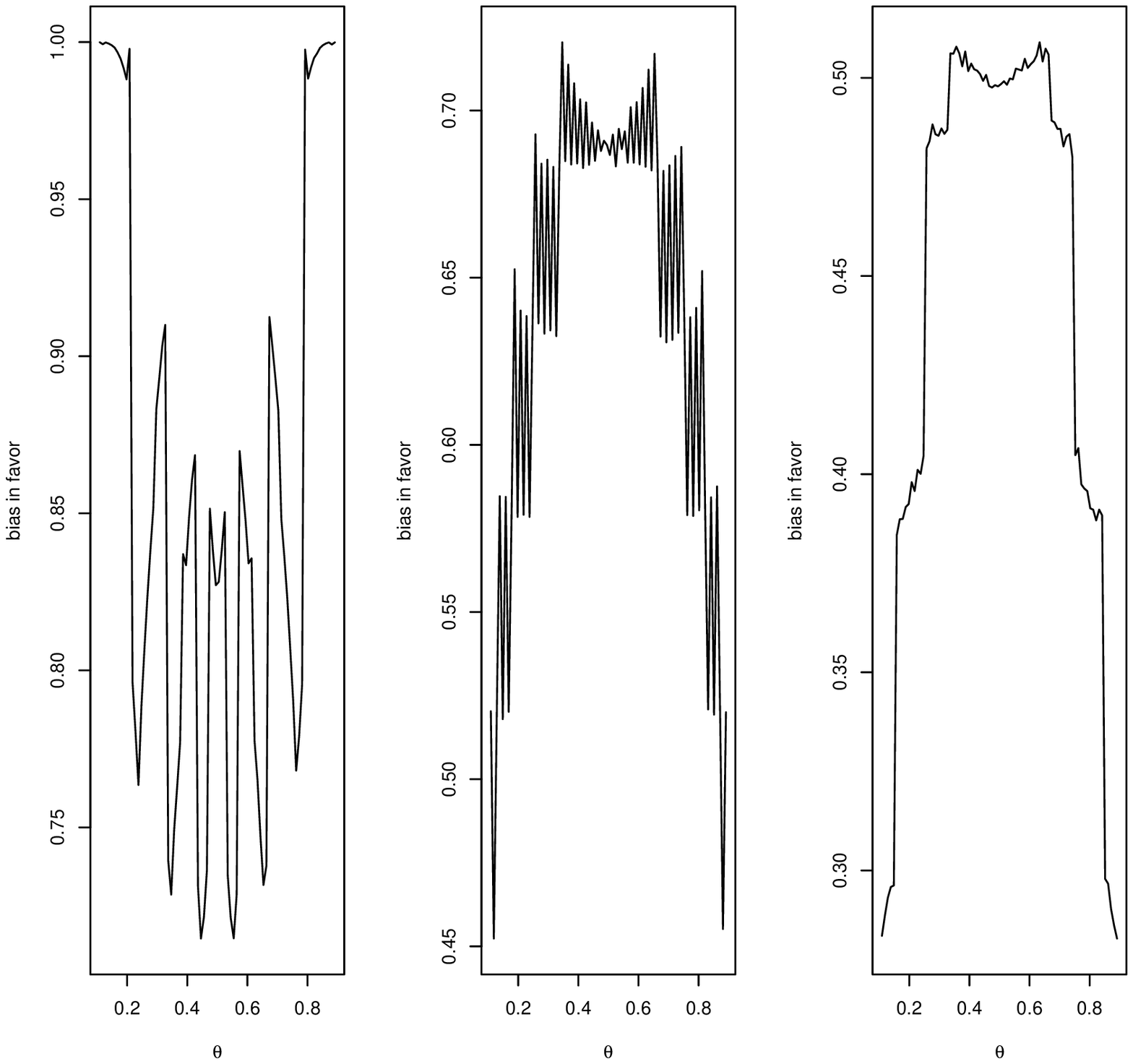}%
\caption{Plots of bias in favor as a function of $\theta$ for $n=10,50$ and
$100$ when using beta$(1,1)$ prior with $\delta=0.1$.}%
\label{binomfig3}%
\end{center}
\end{figure}

\begin{figure}
[p]
\begin{center}
\includegraphics[
height=2.2572in,
width=2.2572in
]%
{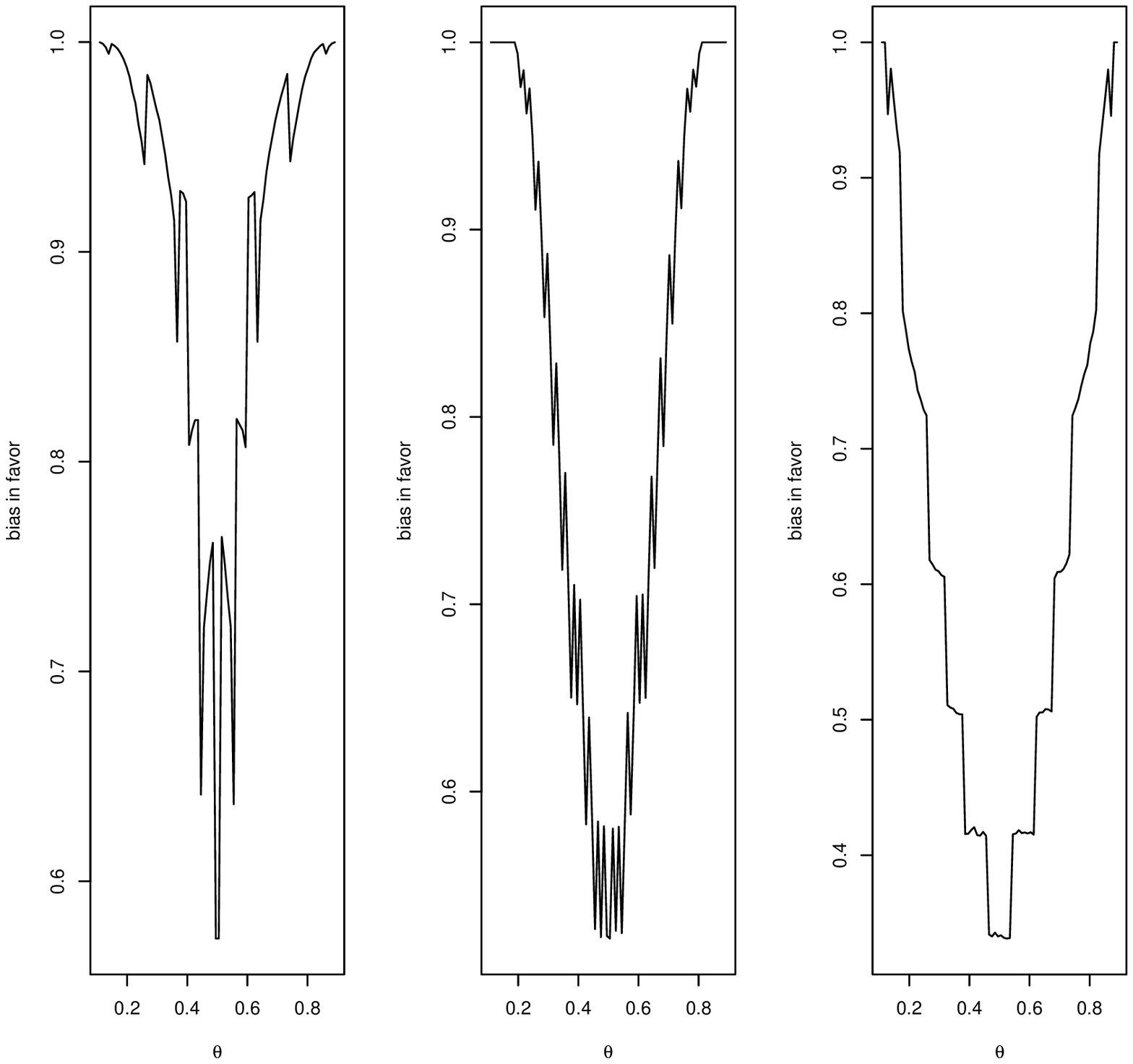}%
\caption{Plots of bias in favor as a function of $\theta$ for $n=10,50$ and
$100$ when using beta$(5,5)$ prior with $\delta=0.1$.}%
\label{binomfig4}
\end{center}
\end{figure}

Now consider bias in favor using (\ref{biasfavest2}). When $(\alpha_{0}%
,\beta_{0})=(1,1)$ and $\delta=0.1,$ Figure \ref{binomfig3} gives the plots of
the bias in favor with $\delta=0.1$ for $n=10$ (max. = $1.00$, average =
$0.84$), $n=50$ (max. = $0.72$, average = $0.51$) and $n=100$ (max. = $0.50$,
average = $0.35$). Therefore, when $n=10$ the maximum probability that $Pl(x)$
contains a false value at least $\delta$ away from the true value is $1,$ when
$n=50$ this probability is $0.72$ and when $n=100$ it is a $0.50.$ When
$(\alpha_{0},\beta_{0})=(5,5)$ Figure \ref{binomfig4} gives the plots of the
bias in favor for $n=10$ (max. = $1.00$, average = $0.68$), for $n=50$ (max. =
$1.00$, average = $0.71$) and for $n=100$ (max. = $1.00$, average = $0.49$).
So in this case the maximum probability that $Pl(x)$ contains a false value at
least $\delta$ away from the true value is always $1,$ but when averaged with
respect to the prior the values are considerably less. It is necessary to
either increase $n$ or $\delta$ to decrease bias in favor. For example, with
$(\alpha_{0},\beta_{0})=(5,5),$ $\delta=0.1$ and $n=400$ the maximum bias in
favor is $0.02$ and the average bias in favor is $0.02$ and when $n=600$ these
quantities equal 0 to two decimals. When $\delta=0.2$ and $n=50$ the maximum
bias in favor is $0.29$ and the average bias in favor is $0.11$ and when
$n=100$ the maximum bias in favor is $0.01$ and the average bias in favor is
$0.01.\smallskip$%
%TCIMACRO{\FRAME{fpFU}{2.2035in}{2.2035in}{0pt}{\Qcb{Plots of bias in favor as
%a function of $\theta$ for $n=10,50$ and $100$ when using beta$(1,1)$ prior
%with $\delta=0.1$.}}{\Qlb{binomfig3}}{binomfig3.eps}%
%width 2.2035in;  height 2.2035in;  depth 0pt;  original-width 6.9998in;
%original-height 6.9998in;  cropleft "0";  croptop "1";  cropright "1";
%filename 'binomialexample/binomfig3.eps';file-properties "XNPEU";}}}%
%BeginExpansion
%EndExpansion%
%TCIMACRO{\FRAME{fpFU}{2.2572in}{2.2572in}{0pt}{\Qcb{Plots of bias in favor as
%a function of $\theta$ for $n=10,50$ and $100$ when using beta$(5,5)$ prior
%with $\delta=0.1$.}}{}{binomfig4.eps}{\special{ language "Scientific Word";
%type "GRAPHIC";  maintain-aspect-ratio TRUE;  display "USEDEF";
%valid_file "F";  width 2.2572in;  height 2.2572in;  depth 0pt;
%original-width 6.9998in;  original-height 6.9998in;  cropleft "0";
%filename 'binomialexample/binomfig4.eps';file-properties "XNPEU";}}}%
%BeginExpansion
%EndExpansion

\noindent\textbf{Example 3.} \textit{Location-scale normal - quantiles.}

Suppose $x=(x_{1},\ldots,x_{n})$ is a sample from $N(\mu,\sigma^{2})$ with
$(\mu,\sigma^{2})\in R^{1}\times(0,\infty)$ unknown with prior $\mu
\,|\,\sigma^{2}\sim N(\mu_{0},\tau_{0}^{2}\sigma^{2}),\sigma^{-2}\sim
\,$gamma$_{\text{rate}}(\alpha_{0},\beta_{0})$. The hyperparameters $(\mu
_{0},\tau_{0}^{2},\alpha_{0},\beta_{0})$ can be obtained via an elicitation
as, for example, discussed in Evans and Tomal (2018) for the more general
regression model. This example is easily generalized to the regression
context. A MSS\ is $T(x)=(\bar{x},||x-\bar{x}1||^{2})$ with the posterior
distribution given by $\mu\,|\,\sigma^{2},T(x)\sim N(\mu_{0x},\left(
n+1/\tau_{0}^{2}\right)  ^{-1}\sigma^{2}),\sigma^{-2}\,|\,T(x)\sim
\,$gamma$_{rate}\left(  \alpha_{0}+n/2,\beta_{0x}\right)  $ where $\mu
_{0x}=(n+1/\tau_{0}^{2})^{-1}(n\bar{x}+\mu_{0}/\tau_{0}^{2})$ and $\beta
_{0x}=\beta_{0}+||x-\bar{x}1||^{2}/2+n(\bar{x}-\mu_{0})^{2}/2(n\tau_{0}%
^{2}+1).$

Suppose interest is in the $\gamma$-th quantile $\psi=\Psi(\mu,\sigma^{2}%
)=\mu+\sigma z_{\gamma},$ where $z_{\gamma}=\Phi^{-1}(\gamma).$ To determine
the bias for or against $\psi$ we need the prior and posterior of $\psi$ which
in this case cannot be worked out in closed form. It is easy, however, to work
with the discretized $\psi$ by simply generating from the prior and posterior
of $(\mu,\sigma^{2}),$ estimate the contents of the relevant intervals and
then approximate the relative belief ratio using these. A natural approach to
the discretization is to base it on the prior mean $E(\psi)=\mu_{0}+\beta
_{0}^{1/2}(\Gamma(\alpha_{0}-1/2)/\Gamma(\alpha_{0}))z_{\gamma}$ and variance
$Var(\psi)=E(\psi^{2})-(E(\psi))^{2}$ where $E(\psi^{2})=(z_{\gamma}^{2}%
+\tau_{0}^{2})\beta_{0}/(\alpha_{0}-1).$ So for a given $\delta,$ we
discretize using $2k+1$ intervals $(E(\psi)+i\delta,E(\psi)+(i+1)\delta]$
where $k=cSD(\psi)/\delta$ and $c$ is chosen so that the collection of
intervals covers the effective support of $\psi$ which is easily assessed as
part of the simulation. For example, with the prior given by hyperparameters
$\mu_{0}=0,\tau_{0}^{2}=1,\alpha_{0}=2,\beta_{0}=1$ and $\gamma=0.5,\delta
=0.1,c=5,$ then $k=50$ and, on generating $10^{5}$ values from the prior,
these intervals contained $99,699$ of the values and with $c=6,$ then $k=60$
and these intervals contained $99,901$ of the generated values. Similar
results are obtained for more extreme quantiles and this is because the
intervals shift with the quantile.

For the bias against for estimation the value of $M(RB_{\Psi}(\psi
\,|\,X)\leq1\,|\,\psi)$ is needed for a range of $\psi$\ values. For this we
need to generate from the conditional prior distribution of $T$ given
$\Psi(\mu,\sigma^{2})=\psi$ and an algorithm for generating from the
conditional prior of $(\mu,\sigma^{2})$ given $\psi$ is needed. Putting
$\nu=1/\sigma^{2},$ the transformation $(\mu,\nu)\rightarrow(\psi,\nu
)=(\mu+\nu^{-1/2}z_{\gamma},\nu)$ has Jacobian equal to 1, so the conditional
prior distribution of $\nu\,|\,\psi$ has density proportional to $\nu
^{\alpha_{0}-1/2}\exp\{-\beta_{0}\nu\}\exp\{-\nu\left(  \psi-\mu_{0}%
-\nu^{-1/2}z_{\gamma}\right)  ^{2}/2\tau_{0}^{2}\}.$ The following gives a
rejection algorithm for generating from this distribution:\smallskip

\noindent1. generate $\nu\sim$ gamma$(\alpha_{0}+1/2,\beta_{0}),$ 2. generate
$u\sim$ unif$(0,1)$ independent of $\nu,$ 3. if $u\leq\exp\{-\nu\left(
\psi-\mu_{0}-\nu^{-1/2}z_{\gamma}\right)  ^{2}/2\tau_{0}^{2}\}$ return $\nu,$
else go to 1.\smallskip

\begin{figure}
[ptb]
\begin{center}
\includegraphics[
height=2.2243in,
width=2.2243in
]%
{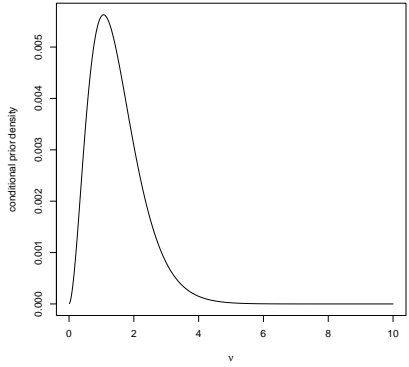}%
\caption{Conditional prior density of $\nu=1/\sigma^{2}$ given $\psi=2$ when
$\gamma=0.95$ and $\mu_{0}=0,\tau_{0}^{2}=1,\alpha_{0}=2,\beta_{0}=1.$}%
\label{condpriorsnu}%
\end{center}
\end{figure}

\noindent As $\psi$ moves away from the prior expected value $E(\psi)$ this
algorithm becomes less efficient but even when the expected number of
iterations is 86 (when $\gamma=0.95,\psi=12),$ generating a sample of $10^{4}$
is almost instantaneous. Figure \ref{condpriorsnu} is a plot of the
conditional prior of $\nu$ given that $\psi=2.$ After generating $\nu$ then
generate $||x-\bar{x}1||^{2}\sim$ $\nu^{-1}$chi-squared$(n-1)$ and $\bar
{x}\sim N(\psi-\nu^{-1/2}z_{\gamma},\nu^{-1}/n)$ to complete the generation of
a value from $M_{T}(\cdot\,|\,\psi).$%
%TCIMACRO{\FRAME{ftbpFU}{2.2243in}{2.2243in}{0pt}{\Qcb{Conditional prior
%density of $\nu=1/\sigma^{2}$ given $\psi=2$ when $\gamma=0.95$ and $\mu
%_{0}=0,\tau_{0}^{2}=1,\alpha_{0}=2,\beta_{0}=1.$}}{\Qlb{condpriorsnu}%
%}{locscalecondprior.eps}{\special{ language "Scientific Word";
%type "GRAPHIC";  maintain-aspect-ratio TRUE;  display "USEDEF";
%valid_file "F";  width 2.2243in;  height 2.2243in;  depth 0pt;
%original-width 6.9998in;  original-height 6.9998in;  cropleft "0";
%filename 'quantileexample/locscalecondprior.eps';file-properties "XNPEU";}}}%
%BeginExpansion
%EndExpansion

\begin{figure}
[ptb]
\begin{center}
\includegraphics[
height=2.5157in,
width=2.5157in
]%
{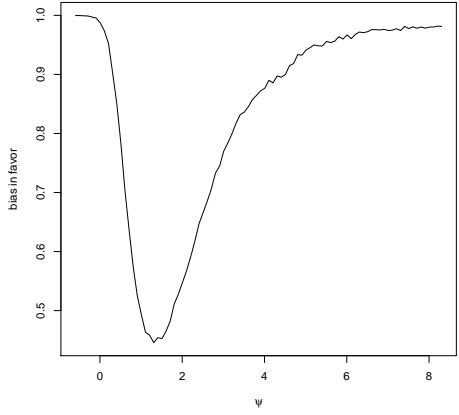}%
\caption{The bias in favor as a function of $\psi$ when $n=10,\delta=0.5$ and
using a prior with hyperparameters $\mu_{0}=0,\tau_{0}^{2}=1,\alpha
_{0}=2,\beta_{0}=1.$}%
\label{biasfavor95}%
\end{center}
\end{figure}

The bias against as a function of $\psi=\mu+\sigma z_{0.95},$ has maximum
value $0.151$ when $n=10$ and so $Pl_{\Psi}(x)$ is a $0.849$-confidence region
for $\psi$ while the average bias against is $0.104$ so the Bayesian coverage
is $0.896.$ Table \ref{coverages} gives the coverages for other values of $n$
as well. Figure \ref{biasfavor95} is a plot of the bias in favor as a function
of $\psi$ with $\delta=\pm0.5$ and $n=10.$ The average bias in favor is
$0.629304.$ When $n=50$ the average bias in favor is $0.3348178.$%

%TCIMACRO{\TeXButton{B}{\begin{table}[tbp] \centering}}%
%BeginExpansion
\begin{table}[tbp] \centering
%EndExpansion%
\begin{tabular}
[c]{|ccc|}\hline
$n$ & \multicolumn{1}{|c}{Frequentist coverage} &
\multicolumn{1}{|c|}{Bayesian coverage}\\\hline
\multicolumn{1}{|c|}{$10$} & \multicolumn{1}{c|}{$0.849$} & $0.896$\\
\multicolumn{1}{|c|}{$20$} & \multicolumn{1}{c|}{$0.895$} & $0.927$\\
\multicolumn{1}{|c|}{$50$} & \multicolumn{1}{c|}{$0.934$} & $0.958$\\
\multicolumn{1}{|c|}{$100$} & \multicolumn{1}{c|}{$0.955$} & $0.973$\\\hline
\end{tabular}
\caption{Coverage probabilities for $Pl_{\psi}(x)$  for the $0.95$ quantile in Example 2.}\label{coverages}%
%TCIMACRO{\TeXButton{E}{\end{table}}}%
%BeginExpansion
\end{table}%
%EndExpansion%
%TCIMACRO{\FRAME{ftbpFU}{2.5157in}{2.5157in}{0pt}{\Qcb{The bias in favor as a
%function of $\psi$ when $n=10,\delta=0.5$ and using a prior with
%hyperparameters $\mu_{0}=0,\tau_{0}^{2}=1,\alpha_{0}=2,\beta_{0}=1.$}%
%}{\Qlb{biasfavor95}}{biasfavor95.eps}{\special{ language "Scientific Word";
%type "GRAPHIC";  maintain-aspect-ratio TRUE;  display "USEDEF";
%valid_file "F";  width 2.5157in;  height 2.5157in;  depth 0pt;
%original-width 6.9998in;  original-height 6.9998in;  cropleft "0";
%filename 'quantileexample/biasfavor95.eps';file-properties "XNPEU";}}}%
%BeginExpansion
%EndExpansion

The case $\gamma=0.50,$ so $\psi=\Psi(\mu,\sigma^{2})=\mu,$ is also of
interest$.$ For $n=10$ then $Pl_{\Psi}\left(  x\right)  $ has $0.878$
frequentist coverage and $0.926$ Bayesian coverage, when $n=20$ the coverages
are $0.916$ and $0.952$ while when $n=50$ the coverages are $0.950$ and
$0.973.$ When $n=10,\delta=0.5$ the average bias in favor is $0.619,$ when
$n=20$ this is $0.4206$ and for $n=100$ the average bias in favor is
$0.091.$\textit{\smallskip}

\noindent\textbf{Example 4.} \textit{Normal Regression - prediction.}

Prediction problems have some unique aspects when compared to inferences about
parameters. To see this consider first the location normal model of Example 1
and suppose the problem is to make inference about a future value $y\sim
N(\mu,\sigma_{0}^{2}).$ The prior predictive distribution is $y\sim N(\mu
_{0},\tau_{0}^{2}+\sigma_{0}^{2})$ and the posterior predictive is $y\sim
N(\mu_{x},\sigma_{n}^{2}+\sigma_{0}^{2})$ where $\mu_{x}=\sigma_{n}^{2}%
(n\bar{x}/\sigma_{0}^{2}+\mu_{0}/\tau_{0}^{2}),\sigma_{n}^{2}=\left(
n/\sigma_{0}^{2}+1/\tau_{0}^{2}\right)  ^{-1}$ and so
\[
RB(y\,|\,\bar{x})=\left(  \frac{\tau_{0}^{2}+\sigma_{0}^{2}}{\sigma_{n}%
^{2}+\sigma_{0}^{2}}\right)  ^{1/2}\exp\left\{  -\frac{1}{2}\left[
\frac{(y-\mu_{x})^{2}}{\sigma_{n}^{2}+\sigma_{0}^{2}}-\frac{(y-\mu_{0})^{2}%
}{\tau_{0}^{2}+\sigma_{0}^{2}}\right]  \right\}  .
\]

For a given $y$ the bias against is given by $M(RB(y\,|\,\bar{x})\leq1\,|\,y)$
and for this we need the conditional prior predictive of $\bar{x}\,|\,y.$ The
joint prior predictive is $(\bar{x},\,y)\sim N_{2}(\mu_{0}1_{2},\Sigma_{0})$
where
\[
\Sigma_{0}=\left(
\begin{array}
[c]{cc}%
\tau_{0}^{2}+\sigma_{0}^{2}/n & \tau_{0}^{2}\\
\tau_{0}^{2} & \tau_{0}^{2}+\sigma_{0}^{2}%
\end{array}
\right)
\]
and so $\bar{x}\,|\,y\sim N(\mu_{0}+\tau_{0}^{2}(y-\mu_{0})/(\tau_{0}%
^{2}+\sigma_{0}^{2}),\sigma_{0}^{2}\left(  \tau_{0}^{2}/(\tau_{0}^{2}%
+\sigma_{0}^{2})+1/n\right)  ).$ From this we see that, as $n\rightarrow
\infty$ the conditional prior distribution of $\mu_{x}\,|\,y$ converges to the
$N\left(  \mu_{0}+\tau_{0}^{2}(y-\mu_{0})/(\tau_{0}^{2}+\sigma_{0}^{2}%
),\sigma_{0}^{2}\tau_{0}^{2}/(\tau_{0}^{2}+\sigma_{0}^{2})\right)  $
distribution. Then with $Z\sim N(0,1)$ and $r=\tau_{0}^{2}/\sigma_{0}^{2}$,
putting $d((y-\mu_{0})/\sigma_{0},r)=\left(  1+1/r\right)  \log\left(
1+r\right)  +r^{-1}(y-\mu_{0})^{2}/\sigma_{0}^{2}$
\[
M(RB(y\,|\,\bar{x})\leq1\,|\,y)\rightarrow1-P\left(  Z\in\left[
\begin{array}
[c]{c}%
r^{-1/2}\left(  1+r\right)  ^{-1/2}\left(  \frac{y-\mu_{0}}{\sigma_{0}%
}\right)  \pm\\
d^{1/2}\left(  \frac{y-\mu_{0}}{\sigma_{0}},r\right)
\end{array}
\right]  \right)
\]
as $n\rightarrow\infty.$ So the bias against does not go to 0 as
$n\rightarrow\infty$ and there is a limiting lower bound to the prior
probability that evidence in favor of a specific $y$ will not be obtained.
This baseline is dependent on both $(y-\mu_{0})/\sigma_{0}$ and $r$. As
$r=\tau_{0}^{2}/\sigma_{0}^{2}\rightarrow\infty$ this baseline bias against
goes to 0 and so it is necessary to ensure that the prior variance is not too
small. Table (\ref{biasagpred}) gives some values for the bias against and it
is seen that if $\tau_{0}^{2}/\sigma_{0}^{2}$ is too small, then there is
substantial bias against even when $y$ is a reasonable value from the
distribution. When $\tau_{0}^{2}/\sigma_{0}^{2}=1,(y-\mu_{0})/\sigma_{0}=0$
and $n=10$ the bias against is computed to be $0.248$ which is quite close to
the baseline so increasing sample size will not reduce bias against by much
and similar results are obtained for the other cases.%

%TCIMACRO{\TeXButton{B}{\begin{table}[tbp] \centering}}%
%BeginExpansion
\begin{table}[tbp] \centering
%EndExpansion
$%
\begin{tabular}
[c]{|l|l|l|}\hline
$\tau_{0}^{2}/\sigma_{0}^{2}$ & bias against $\frac{y-\mu_{0}}{\sigma_{0}}=0$
& bias against $\frac{y-\mu_{0}}{\sigma_{0}}=1$\\\hline
\multicolumn{1}{|r|}{$1$} & \multicolumn{1}{|c|}{$0.239$} &
\multicolumn{1}{|c|}{$0.213$}\\
\multicolumn{1}{|r|}{$10$} & \multicolumn{1}{|c|}{$0.104$} &
\multicolumn{1}{|c|}{$0.100$}\\
\multicolumn{1}{|r|}{$100$} & \multicolumn{1}{|c|}{$0.031$} &
\multicolumn{1}{|c|}{$0.031$}\\
\multicolumn{1}{|r|}{$1/2$} & \multicolumn{1}{|c|}{$0.270$} &
\multicolumn{1}{|c|}{$0.263$}\\
\multicolumn{1}{|r|}{$1/100$} & \multicolumn{1}{|c|}{$0.316$} &
\multicolumn{1}{|c|}{$0.460$}\\\hline
\end{tabular}
$
\caption{Baseline bias against values for prediiction for location normal in Example 4.}\label{biasagpred}%
%TCIMACRO{\TeXButton{E}{\end{table}}}%
%BeginExpansion
\end{table}%
%EndExpansion

\begin{figure}
[t]
\begin{center}
\includegraphics[
height=2.4837in,
width=2.4837in
]%
{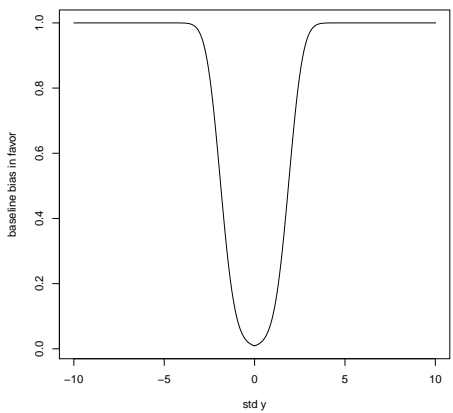}%
\caption{Plot of the baseline bias in favor for values of $(y-\mu_{0}%
)/\sigma_{0}$ when $\tau_{0}^{2}/\sigma_{0}^{2}=1$ when $\delta=5.$}%
\label{biasforpred}%
\end{center}
\end{figure}

Now consider bias in favor of $y$, namely, $M(RB(y\,|\,\bar{x})\geq
1\,|\,y\pm\delta)$ for some choice of $\delta.$ False values for $y$
correspond to values in the tails so we consider, for example, $y+\delta$ as a
value in the central region of the prior and then a large value of $\delta$
puts $y$ in the tails. Again the bias in favor has a baseline value as
$n\rightarrow\infty.$ A similar argument leads to the bias in favor of $y$
satisfying%
\begin{align*}
&  M(RB(y\,|\,\bar{x})\geq1\,|\,y\pm\delta)\\
&  \rightarrow P\left(  Z\in\left[  r^{-1/2}\left(  1+r\right)  ^{-1/2}\left(
\frac{y-\mu_{0}}{\sigma_{0}}\pm r\frac{\delta}{\sigma_{0}}\right)  \pm
d^{1/2}\left(  \frac{y-\mu_{0}}{\sigma_{0}},r\right)  \right]  \right)  .
\end{align*}
Figure \ref{biasforpred} is a plot of $\sup M(RB(y\,|\,\bar{x}%
)\geq1\,|\,y\pm\delta).$ So the bias in favor is reasonably low for central
values of $y$ but it is to be noted that once again there is a trade-off as
when $\tau$ increases the bias in favor goes to 1.%
%TCIMACRO{\FRAME{ftFU}{2.4837in}{2.4837in}{0pt}{\Qcb{Plot of the baseline bias
%in favor for values of $(y-\mu_{0})/\sigma_{0}$ when $\tau_{0}^{2}/\sigma
%_{0}^{2}=1$ when $\delta=5.$}}{\Qlb{biasforpred}}{fig13.eps}%
%width 2.4837in;  height 2.4837in;  depth 0pt;  original-width 6.9998in;
%original-height 6.9998in;  cropleft "0";  croptop "1";  cropright "1";
%filename 'regressionexample/fig13.eps';file-properties "XNPEU";}}}%
%BeginExpansion
%EndExpansion

Prediction plays a bigger role in regression problems but we can expect the
same issues to apply as in the location problem. Suppose $y\sim N_{n}%
(X\beta,\sigma^{2}I)$ where $X\in R^{n\times k}$ is of rank $k,$
$(\beta,\sigma^{2})\in R^{k}\times(0,\infty)$ is unknown, our interest is in
predicting a future value $y_{new}\sim N(w^{t}\beta,\sigma^{2})$ for some
fixed known $w$ and, putting $\nu=1/\sigma^{2},$ the conjugate prior
$\beta\,|\,\nu\sim N_{k}(\beta_{0},\nu^{-1}\Sigma_{0})\,,\nu\sim
\,$gamma$_{\text{rate}}(\alpha_{0},\eta_{0})$ is used. Specifying the
hyperparameters $(\beta_{0},\Sigma_{0},\alpha_{0},\eta_{0})$ can be carried
out using an elicitation algorithm such as that discussed in Evans and Tomal (2018).

For the bias calculations it is necessary to generate values of the MSS
$(b,s^{2})=((X^{t}X)^{-1}X^{t}y,||y-Xb||^{2})$ from the conditional prior
predictive $M(\cdot\,|\,y_{new}).$ This is accomplished by generating from the
conditional prior of $(\beta,\nu)\,|\,y_{new}$ and then generating $b\sim
N_{k}(\beta,\nu^{-1}(X^{t}X)^{-1})$ independent of $s^{2}\sim\nu^{-1}%
\,$chi-squared$(n-k).$ The conditional prior of $(\beta,\nu)\,|\,y_{new}$ is
proportional to
\begin{align*}
&  \nu^{\alpha_{0}-1/2}\exp\{-\eta_{0}(y_{new})\nu\}\times\\
&  \nu^{k/2}\exp\left\{  -\frac{\nu}{2}\left(  \beta-\left(  \Sigma_{0}%
^{-1}+ww^{t}\right)  ^{-1}(\Sigma_{0}^{-1}\beta_{0}+y_{new}w)\right)
^{t}\left(  \Sigma_{0}^{-1}+ww^{t}\right)  \left(  \cdot\right)  \right\}
\end{align*}
where, using $\left(  \Sigma_{0}^{-1}+ww^{t}\right)  ^{-1}=\Sigma_{0}%
-(1+w^{t}\Sigma_{0}w)^{-1}\Sigma_{0}ww^{t}\Sigma_{0},\eta_{0}(y_{new}%
)=\eta_{0}+(1+w^{t}\Sigma_{0}w)^{-1}(w^{t}\beta-y_{new})^{2}/2.$ So generating
$(\beta,\nu)\,|\,y_{new}$ is accomplished via $\nu\sim\,$gamma$_{\text{rate}%
}(\alpha_{0}+1/2,\eta_{0}(y_{new})),$
\[
\beta\,|\,\nu\sim N_{k}\left(  \left(  I-\frac{\Sigma_{0}ww^{t}}{1+w^{t}%
\Sigma_{0}w}\right)  (\beta_{0}+y_{new}\Sigma_{0}w),\nu^{-1}\left(  \Sigma
_{0}-\frac{\Sigma_{0}ww^{t}\Sigma_{0}}{1+w^{t}\Sigma_{0}w}\right)  \right)  .
\]
For each generated $(b,s^{2})$ it is necessary to compute the relative belief
ratio $RB(y_{new}\,|\,b,s^{2})$ and determine if it is less than or equal to
$1.$ There are closed forms for the prior and conditional densities of
$y_{new}$ since
\begin{align*}
y_{new}  &  \sim w^{t}\beta_{0}+\left\{  \eta_{0}(1+w^{t}\Sigma_{0}%
w)/\alpha_{0}\right\}  ^{1/2}t_{2\alpha_{0}},\\
y_{new}\,|\,(b,s^{2})  &  \sim w^{t}\beta_{0}(b,s^{2})+\left\{  \frac{\eta
_{0}(b,s^{2})(1+w^{t}\left(  \Sigma_{0}^{-1}+X^{t}X\right)  ^{-1}w)}%
{\alpha_{0}+n/2}\right\}  ^{1/2}t_{2\alpha_{0}+n}%
\end{align*}
where $t_{\lambda}$ denotes a Student$(\lambda)$ random variable and%
\begin{align*}
\beta_{0}(b,s^{2})  &  =\left(  \Sigma_{0}^{-1}+X^{t}X\right)  ^{-1}\left(
\Sigma_{0}^{-1}\beta_{0}+X^{t}Xb\right) \\
\eta_{0}(b,s^{2})  &  =\eta_{0}+\left[  s^{2}+||Xb||^{2}+||\Sigma_{0}%
^{-1}\beta_{0}||^{2}-\beta_{0}(b,s^{2})^{t}\left(  \Sigma_{0}^{-1}%
+X^{t}X\right)  \beta_{0}(b,s^{2})\right]  /2.
\end{align*}
These results permit the calculation of the biases as in the location problem.

\section{Conclusions}

There are several conclusions that can be drawn from the discussion here.
First, it is necessary to take bias into account when considering Bayesian
procedures and currently this is generally not being done. Depending on the
purpose of the study, some values concerning both bias against and bias in
favor need to be quoted as these are figures \ of merit for the study. The
approach to Bayesian inferences via a characterization of evidence makes this
relatively straight-forward conceptually. Second, frequentism plays a role in
Bayesian statistical reasoning, not through the inferences, but rather through
the design as it how we determine and control the biases. Overall this makes
sense because, before the data is seen, it is natural to be concerned about
what inferences can be reliably drawn. Once the data is observed, however, it
is the evidence in this data set that matters and not the evidence in the data
sets not seen. Still, if we ignore the latter it may be that the existence of
bias makes the inferences drawn of very low quality. Third, the results
concerning the standard p-value in Example 1 can be seen to apply quite
generally and this makes any discussion about how to characterize and measure
evidence of considerable importance. The principle of evidence makes a
substantial contribution in this regard as was shown in a variety of results.
The major purpose of this paper, however, is to deal with a key criticism of
Bayesian methodology, namely, that inferences can be biased because of their
dependence on the subjective beliefs of the analyst. This criticism is
accepted, but we also assert that this can be dealt with in a logical and
scientific fashion as has been demonstrated in this paper.

\section{References}

\noindent Baskurt, Z. and Evans, M. (2013) Hypothesis assessment and
inequalities for Bayes factors and relative belief ratios. Bayesian Analysis,
8, 3, 569-590.

\noindent Berger, J.O. and Selke, T. (1987) Testing a point null hypothesis:
the irreconcilability of p values and evidence. Journal of the American
Statistical Association, 82, 397, 112-122.

\noindent Berger, J.O. and Delampady, M. (1987) Testing precise hypotheses.
Statistical Science, 2, 3, 317-335.

\noindent Cousins, R.D. (2017) The Jeffreys--Lindley paradox and discovery
criteria in high energy physics. Synthese, 194, 2, 395--432.

\noindent Evans, M. (2015) Measuring Statistical Evidence Using Relative
Belief. Monographs on Statistics and Applied Probability 144, CRC Press.

\noindent Evans, M., Guttman, I. and Li, P. (2017) Prior elicitation,
assessment and inference with a Dirichlet prior. Entropy 2017, 19(10), 564; doi:10.3390/e1910056.

\noindent Evans, M. and Tomal, J. (2018) Multiple testing via relative belief
ratios. FACETS, 3: 563-583, DOI: 10.1139/facets-2017-0121.

\noindent Gu, Y., Li, W. Evans, M. and Englert, B-G. (2019) Very strong
evidence in favor of quantum mechanics and against local hidden variables from
a Bayesian analysis. Physical Review A 99, 022112(1-17).

\noindent Robert, C. P. (2014) On the Jeffreys-Lindley paradox. Philosophy of
Science, 81, 216--232.

\noindent Shafer, G. (1982)\emph{ }Lindley's paradox (with discussion).
Journal of the American Statistical Association, 77, 378, 325-351.

\noindent Spanos, A. Who should be afraid of the Jeffreys-Lindley paradox?
Philosophy of Science, 80, 1, 73 - 93.

\noindent Sprenger, J. (2013) Testing a precise null hypothesis: The case of
Lindley's paradox. Philosophy of Science, 80, 733-744.

\noindent Villa, C. and Walker, S. (2017) On the mathematics of the
Jeffreys--Lindley paradox. Communications in Statistics - Theory and Methods,
46, 24, 12290-12298.

\section{Appendix}

\noindent\textbf{Proof of Theorem 1. }The Savage-Dickey ratio result implies
$RB_{\Psi}(\psi_{\ast}\,|\,x)=m_{\psi_{\ast}}(x)/m(x)$ and note $R(\psi_{\ast
})=\{x:m_{\psi_{\ast}}(x)\leq m(x)\}.$ Now $\mathcal{X}_{1}=\{x:I_{R(\psi
_{\ast})}(x)-I_{D(\psi_{\ast})}(x)<0\}=\{x:I_{R(\psi_{\ast})}-I_{D(\psi_{\ast
})}(x)<0,m_{\psi_{\ast}}(x)>m(x)\}$ and $\mathcal{X}_{2}=\{x:I_{R(\psi_{\ast
})}(x)-I_{D(\psi_{\ast})}(x)>0\}=\{x:I_{R(\psi_{\ast})}(x)-I_{D(\psi_{\ast}%
)}(x)\geq0,m_{\psi_{\ast}}(x)\leq m(x)\}.$ Then $M(R(\psi_{\ast}%
))-M(D(\psi_{\ast}))=\int_{\mathcal{X}_{1}}(I_{R(\psi_{\ast})}(x)-I_{D(\psi
_{\ast})}(x))\,$\newline$M(dx)+\int_{\mathcal{X}_{2}}(I_{R(\psi_{\ast}%
)}(x)-I_{D(\psi_{\ast})}(x))\,M(dx)\geq M(R(\psi_{\ast})\,|\,\psi_{\ast
})-M(D(\psi_{\ast})\,|\,\psi_{\ast})\geq0$ establishing (i). Also,
$M(D(\psi_{\ast}))=M(D(\psi_{\ast})\,|\,\psi_{\ast})\Pi_{\Psi}(\{\psi_{\ast
}\})+\int_{\Psi\backslash\{\psi_{\ast}\}}M(D(\psi_{\ast})\,|\,$\newline%
$\psi)\,\Pi_{\Psi}(d\psi)$ and the integral is the prior probability of not
getting evidence in favor of $\psi_{\ast}$ when it is false and this
establishes (ii).\smallskip

\noindent\textbf{Proof of Theorem 2.} Now $E_{\Pi_{\Psi}}(M(\psi_{\ast}\notin
C(X)))=E_{\Pi_{\Psi}^{2}}\left(  M(\psi_{\ast}\notin C(X)\,|\,\psi)\right)
=E_{\Pi_{\Psi}^{2}}\left(  M(D(\psi_{\ast}))\,|\,\psi)\right)  =\int_{\Psi
}M(D(\psi_{\ast}))\,\Pi_{\Psi}(d\psi_{\ast})$ and (i) follows from Theorem 1.
Also, $\int_{\Psi}M(D(\psi_{\ast}))\,\Pi_{\Psi}(d\psi_{\ast})=E_{\Pi_{\Psi}%
}(\int_{\Psi}M(D(\psi_{\ast})\,|\,\psi)\,\Pi_{\Psi}(d\psi))=E_{\Pi_{\Psi}%
}(M(D(\psi_{\ast})\,$\newline$|\,\psi_{\ast})\Pi_{\Psi}(\{\psi_{\ast
}\}))+E_{\Pi_{\Psi}}(\int_{\Psi\backslash\{\psi_{\ast}\}}M(D(\psi_{\ast
})\,|\,\psi)\,\Pi_{\Psi}(d\psi))=E_{\Pi_{\Psi}}(M(\psi_{\ast}\notin
C(X)\,|\,\psi_{\ast})$\newline$\Pi_{\Psi}(\{\psi_{\ast}\}))+E_{\Pi_{\Psi}%
}(\int_{\Psi\backslash\{\psi_{\ast}\}}M(\psi_{\ast}\notin C(X)\,\,|\,\psi
)\,\Pi_{\Psi}(d\psi))$ establishing (ii).\smallskip

\noindent\textbf{Proof of Theorem 3.} Now $M(R(\psi_{\ast})\,|\,\psi_{\ast
})=\int I_{R(\psi_{\ast})}(x)\,M_{\psi_{\ast}}(dx)\leq\int I_{R(\psi_{\ast}%
)}(x)\,$\newline$M(dx)=M(R(\psi_{\ast}))=\int_{\Psi}M(R(\psi_{\ast}%
)\,|\,\psi)\,\Pi(d\psi)=M(R(\psi_{\ast})\,|\,\psi_{\ast})\Pi_{\Psi}%
(\{\psi_{\ast}\})+\int_{\Psi\backslash\{\psi_{\ast}\}}M(R(\psi_{\ast
})\,|\,\psi)\,\Pi_{\Psi}(d\psi)$ and so $\Pi_{\Psi}(\{\psi_{\ast}%
\}^{c})M(R(\psi_{\ast})\,|\,\psi_{\ast})\leq\int_{\Psi\backslash\{\psi_{\ast
}\}}$\newline$M(R(\psi_{\ast})\,|\,\psi)\,\Pi_{\Psi}(d\psi)$ which implies
(i). Furthermore, (ii) is implied by
\begin{align*}
&  E_{\Pi_{\Psi}}\left(  M(\psi_{\ast}\notin Pl_{\Psi}(X)\,|\,\psi_{\ast
})\right)  =E_{\Pi_{\Psi}}\left(  M(R(\psi_{\ast})\,|\,\psi_{\ast})\right) \\
&  \leq E_{\Pi_{\Psi}}\left(  \int_{\Psi\backslash\{\psi_{\ast}\}}%
M(R(\psi_{\ast})\,|\,\psi)\,\Pi_{\Psi}(d\psi)/\Pi_{\Psi}(\{\psi_{\ast}%
\}^{c}\right) \\
&  =E_{\Pi_{\Psi}}\left(  \int_{\Psi\backslash\{\psi_{\ast}\}}M(\psi_{\ast
}\notin Pl_{\Psi}(X)\,|\,\psi)\,\Pi_{\Psi}(d\psi)/\Pi_{\Psi}(\{\psi_{\ast
}\}^{c}\right)  .
\end{align*}

\noindent\textbf{Proof of Theorem 4.} It is easy to see that the proof of
Theorem 1 can be modified to show that among all regions $D^{int}(\psi_{\ast
})\subset\mathcal{X}$ satisfying $M(D^{int}(\psi_{\ast})\,|\,\psi_{\ast})\leq
M(RB_{\Psi}(\psi_{\ast}\,|\,X)<1\,|\,\psi_{\ast})$ the prior probability
$M(D^{int}(\psi_{\ast}))$ is maximized by $D^{int}(\psi_{\ast})=\{x:RB_{\Psi
}(\psi_{\ast}\,|\,x)<1\}.$ This clearly implies (i) and (ii) follows
similarly.\smallskip

\noindent\textbf{Proof of Theorem 5.} Now $E_{\Pi_{\Psi}}(M(\psi_{\ast}\in
C(X)))=E_{\Pi_{\Psi}^{2}}\left(  M(\psi_{\ast}\in C(X)\,|\,\psi)\right)
=E_{\Pi_{\Psi}^{2}}\left(  M(D^{c}(\psi_{\ast}))\,|\,\psi)\right)
=E_{\Pi_{\Psi}}(M(D^{c}(\psi_{\ast}))$ and (i) follows from Theorem 1(i).
Also, (ii) is implied by $E_{\Pi_{\Psi}}(M(D^{c}(\psi_{\ast}))=\int_{\Psi
}M(D^{c}(\psi_{\ast})\,|\,\psi_{\ast})\Pi_{\Psi}(\{\psi_{\ast}\})\,\Pi_{\Psi
}(d\psi_{\ast})+\int_{\Psi}\int_{\Psi\backslash\{\psi_{\ast}\}}M(D^{c}%
(\psi_{\ast})\,|\,\psi)\,\Pi_{\Psi}(d\psi)\,\Pi_{\Psi}(d\psi_{\ast}%
)=\int_{\Psi}M(\psi_{\ast}\in C(X)\,|\,\psi_{\ast})\Pi_{\Psi}(\{\psi_{\ast
}\})\,\newline\Pi_{\Psi}(d\psi_{\ast})+\int_{\Psi}\int_{\Psi\backslash
\{\psi_{\ast}\}}M(\psi_{\ast}\in C(X)\,\,|\,\psi)\,\Pi_{\Psi}(d\psi
)\,\Pi_{\Psi}(d\psi_{\ast}).$

\end{document}